\documentclass[11pt,reqno]{amsart}

\usepackage{booktabs}                                 
\usepackage[english]{babel}
\usepackage{latexsym}
\usepackage[utf8]{inputenc}
\usepackage{babel}
\usepackage{fancyhdr}
\usepackage{longtable}
\usepackage{mathrsfs}
\usepackage{geometry}
\usepackage{comment}
\usepackage{textcomp}
\usepackage{caption}
\usepackage{lmodern}
\usepackage{mathrsfs}
 \usepackage[T1]{fontenc}

 \usepackage{multicol}

\usepackage[all,cmtip]{xy}
\usepackage{epsfig}

\usepackage{amsmath,amsfonts,amssymb,amsthm}
\usepackage{graphics}
\usepackage{graphicx}
\usepackage{verbatim}
\usepackage{dsfont}
\usepackage{enumerate}  
\usepackage{pdflscape}
\usepackage{longtable}
\usepackage{booktabs}
\usepackage{colortbl}
\usepackage{systeme} 

\usepackage{tikz-cd}
\usepackage{ytableau}

\usepackage[shortlabels]{enumitem}
\usepackage[breaklinks]{hyperref} 
\hypersetup{
	colorlinks,
	linkcolor={red},
	citecolor={blue},
	urlcolor={red},
}
\usepackage{stmaryrd}
\usepackage{multirow}
\usepackage{longtable}

\usepackage{adjustbox}

\usepackage{xcolor,colortbl}
\definecolor{lavender}{rgb}{0.9, 0.9, 0.98}

\usepackage[]{todonotes}

\newcommand{\red}[1]{{\color{red}{#1}}}

\newcommand{\Z}{\mathbb{Z}}

\newcommand{\mQ}{\mathcal{Q}}

\newcommand{\mF}{\mathcal{F}}

\newcommand{\mE}{\mathcal{E}}
\newcommand{\mU}{\mathcal{U}}

\newcommand{\mZ}{\mathscr{Z}}
\newcommand{\of}{\mathcal{O}}

\newcommand{\Schur}{\Sigma}

\newcommand{\hk}{hyperk\"ahler }

\usepackage{faktor}

\newcommand{\ch}{\mathsf{ch}}

\DeclareMathOperator{\Ext}{Ext}

\DeclareMathOperator{\End}{End}

\DeclareMathOperator{\Sym}{Sym}

\DeclareMathOperator{\Gr}{Gr}

\newtheorem{thm}{Theorem}[section]
\newtheorem{corollary}[thm]{Corollary}

\newtheorem{lemma}[thm]{Lemma}
\newtheorem{proposition}[thm]{Proposition}

\newtheorem*{aim*}{Aim of this paper}

\theoremstyle{definition}
\newtheorem{definition}[thm]{Definition}
\newtheorem{rmk}[thm]{Remark}
\newtheorem{ex}[thm]{Example}

\newtheorem{notz}[thm]{Notation}
\newtheorem{alg}[thm]{Algorithm}

\usepackage{thmtools}
\declaretheoremstyle[
spaceabove=1.5ex, spacebelow=1.5ex,
headfont=\bf,
notefont=\mdseries, notebraces={(}{)},
bodyfont=\normalfont,
headpunct=.,
numberwithin=,
postheadhook=\leavevmode%
qed={$\bullet$}
]{mystyle}

\usepackage[capitalise,nameinlink]{cleveref}
\crefname{table}{Table}{Tables}
\crefname{section}{Section}{Sections}
\crefname{mystyle}{Fano}{Fanos}
\crefname{claim}{Claim}{Claims}
\crefname{rmk}{Remark}{Remarks}
\crefname{workhyp}{WH}{WH}
\crefname{thm}{Theorem}{Theorems}
\crefname{proposition}{Proposition}{Propositions}
\crefname{app}{Appendix}{Appendices}
\crefname{eq}{Equation}{Equations}
\crefname{lemma}{Lemma}{Lemmas}
\crefname{alg}{Algorithm}{Algorithms}
\crefname{ex}{Example}{Examples}
\crefname{conj}{Conjecture}{Conjectures}
\crefname{notz}{Notation}{Notations}
\makeatletter    
\def\l@subsection{\@tocline{1}{0,2pt}{2pc}{8mm}{\ \ }} 
\makeatletter    
\def\l@section{\@tocline{1}{0,2pt}{2pc}{8mm}{\ \ }} 

\newcommand{\crefpart}[2]{%
  \hyperref[#2]{\namecref{#1}~\labelcref*{#1}~(\ref*{#2})}%
}
\author{Alessandro Frassineti}
\address{Dipartimento di Matematica \\
Dipartimento di Eccellenza 2023-2027\\
Universit\`a di Genova\\
Via Dodecaneso 35\\
16146 Genova, Italy}
\email[A.~Frassineti]{alessandro.frassineti@edu.unige.it}

\author{Federico Tufo}
\address{Dipartimento di Matematica \\
Universit\`a di Bologna\\
Piazza di Porta San Donato 5\\
40127 Bologna, Italy}
\email[F.~Tufo]{federico.tufo2@unibo.it}
\usepackage[symbol]{footmisc}

\usepackage{libertine}
\usepackage{libertinust1math} 

\usepackage{float}
\floatstyle{plaintop}
\restylefloat{table}

\usepackage{ytableau}

\title[Modular vector bundles on Hyperk\"ahler manifolds of Debarre--Voisin type]{Modular vector bundles on Hyperk\"ahler manifolds of Debarre--Voisin type}

\begin{document}
\begin{abstract}
    Let $X$ be a very general Debarre--Voisin fourfold. In this article, we prove that all the Schur functors of the restriction of the quotient bundle of $\mathrm{Gr}(6,10)$ to $X$ are modular and polystable vector bundles. We also show that the atomic bundles among them are exactly the symmetric powers of the restriction of the quotient bundle. Moreover, we compute the $\mathrm{Ext}$-groups of different modular vector bundles on $X$, and we find examples with non-trivial dimensional first $\mathrm{Ext}$-group.   
\end{abstract}
\maketitle
\section{Introduction}
Hyperk\"ahler manifolds, or irreducible holomorphic symplectic manifolds, are among the central objects of study in modern algebraic geometry. 
Only a few examples of such varieties are currently known. In every even complex dimension, two deformation families appear: the Hilbert scheme of $n$ points on a $\textrm{K3}$ surface, known as $\textrm{K3}^{[n]}$-type, and the generalized Kummer varieties, known as $\mathrm{Kum}_n$-type. In addition, there exist two sporadic examples in dimensions ten and six, denoted by $\mathrm{OG}10$ and $\mathrm{OG}6$, which arise as desingularizations of certain moduli spaces of semistable sheaves on $\textrm{K3}$ and abelian surfaces, respectively. A posteriori, all these examples admit an interpretation as moduli spaces of sheaves. This naturally leads to the question of whether one can systematically construct families of \hk manifolds as moduli spaces of suitably chosen stable sheaves.

Motivated by this perspective, the study of sheaves on higher-dimensional \hk manifolds has recently attracted increasing attention. A natural starting point is to generalize the rich theory of moduli spaces of sheaves on $\mathrm{K3}$ surfaces, developed since Mukai's seminal work \cite{Muk87}. A key feature that one would like to preserve in this higher-dimensional setting is the good behavior of such moduli spaces in families. Achieving this appears to require restricting to sheaves with special geometric properties. To date, several related notions have been introduced in this direction.

The first approach is due to Verbitsky \cite{Ver96}, who introduced the notion of projectively hyperholomorphic vector bundles on hyperk\"ahler manifolds in order to control deformations along twistor families. Using this framework, he proved that the irreducible component of the moduli space of simple sheaves containing a fixed hyperholomorphic vector bundle naturally carries a hyperk\"ahler structure.

Subsequently, O'Grady \cite{OG22} introduced the notion of modular sheaves. For such sheaves, the wall-and-chamber structure governing stability in the ample (and later K\"ahler) cone behaves analogously to the case of $\textrm{K3}$ surfaces. Later, Markman \cite{Mar24} studied another class of coherent sheaves, subsequently termed atomic by Beckmann \cite{Bec25}. These are characterized by requiring their Mukai vector to lie in the extended Mukai lattice, a condition which is automatically satisfied in the $\textrm{K3}$ case.

These notions are closely related. Any atomic torsion-free sheaf is modular (see \cite[Theorem 1.2]{Mar24}), and any atomic polystable vector bundle is projectively hyperholomorphic (see \cite[Subsec.\ 5.2]{Bec25}). The converse does not hold in general: for instance, the tangent bundle is always modular (\cite{OG22}) but is not atomic for known families of \hk manifolds of dimension greater than two (\cite[Proposition 8.3]{Bec25}). In the special case of $\textrm{K3}^{[2]}$-type manifolds, one can say more. Combining results from \cite{Ver96}, any polystable modular vector bundle is projectively hyperholomorphic.

In this paper, we focus on modular sheaves. For simplicity, we restrict to the case of $\textrm{K3}^{[2]}$-type manifolds, where the definition takes a particularly simple form.

\begin{definition}
Let $X$ be a \hk manifold of $\textrm{K3}^{[2]}$-type. A torsion-free sheaf $\mathcal{E}$ on $X$ is said to be \emph{modular} if there exists a rational number $d$ such that
\[
\Delta(\mathcal{E}) = d \mathsf{c}_2(X),
\]
where
\[
\Delta(\mathcal{E}) := 2r(\mathcal{E})\mathsf{c}_2(\mathcal{E}) - (r(\mathcal{E})-1)\mathsf{c}_1(\mathcal{E})^2 = \ch_1(\mathcal{E})^2 - 2 \ch_0(\mathcal{E})\ch_2(\mathcal{E})
\]
is called the \emph{discriminant} of $\mathcal{E}$.
\end{definition}

As soon as the notion of modular sheaves was introduced in \cite{OG22}, examples were constructed. In particular, O'Grady provided a general method to produce modular sheaves on \hk manifolds of $\textrm{K3}^{[2]}$-type. However, all examples in \textit{loc.\ cit.} are rigid, meaning that their moduli space consists of a single point, and thus they are not suitable for producing new \hk manifolds.

The first examples of non-rigid modular sheaves were constructed by Markman in \cite{Mar24}. Moreover, Bottini \cite{Bot24} constructed a modular sheaf whose moduli space is associated with a \hk manifold of $\mathrm{OG}10$-type. Later, Fatighenti \cite{Fat24} produced further examples using different techniques, which were subsequently generalized by Fatighenti and Onorati \cite{FatOno}. They worked on an explicit model of a \hk manifold of $\textrm{K3}^{[2]}$-type, namely the Beauville--Donagi fourfold, and obtained the following result.

\begin{ex}[Fatighenti--Onorati]
Let $X_{BD}\subseteq \mathrm{Gr}(2,6)$ be the Fano variety of lines of a very general cubic fourfold, which is a \hk manifold of $\textrm{K3}^{[2]}$-type. Let $\widetilde{\mQ}$ be the tautological quotient bundle on $\mathrm{Gr}(2,6)$ and let $\mQ$ be its restriction to $X_{BD}$. Then, for any partition $\lambda = (\lambda_1, \lambda_2, \lambda_3, \lambda_4)$, the associated Schur functor $\Sigma_\lambda \mQ$ is a modular vector bundle. Moreover, it is always polystable and, for small values of $\lambda_i$, it is slope-stable with
\[
\mathrm{ext}^1(\Sigma_\lambda \mQ,\Sigma_\lambda \mQ) \in \{0,20,40\}.
\]
\end{ex}

The Fano variety of lines of a cubic fourfold is one of the few examples of \hk manifolds that can be realized as the zero locus of a section of a homogeneous vector bundle inside a Grassmannian. It is therefore natural to ask whether similar constructions can be carried out in other cases. In this paper, we address this question for the Debarre--Voisin fourfold. In particular, since both the Beauville--Donagi and Debarre--Voisin fourfolds belong to locally complete families of polarized \hk manifolds of $\textrm{K3}^{[2]}$-type, we study modular sheaves on a very general Debarre--Voisin fourfold $X_{DV}\subseteq \mathrm{Gr}(6,10)$, i.e., a smooth Debarre--Voisin fourfold with Picard rank $1$. Our results can be summarized as follows:

\begin{thm}\label[thm]{MainTheorem} [\cref{SchurIsModular}, \cref{rmk:stab}, \cref{prop:Tab}]
Let $X_{DV} \subseteq \mathrm{Gr}(6,10)$ be a very general Debarre--Voisin \hk manifold. Let $\widetilde{\mQ}$ be the tautological quotient bundle on $\mathrm{Gr}(6,10)$ and let $\mQ$ be its restriction to $X_{DV}$. Let $\lambda = (\lambda_1, \lambda_2, \lambda_3, \lambda_4)$ be a partition and consider the associated Schur functor $\Sigma_\lambda \mQ$. Then:
\begin{enumerate}
\item[(1)] $\Sigma_\lambda \mQ$ is modular for every $\lambda$;
\item[(2)] $\Sigma_\lambda \mQ$ is polystable for every $\lambda$;
\item[(3)] the dimensions of $\Ext^*(\Sigma_\lambda \mQ,\Sigma_\lambda \mQ)$ are given in \cref{MainTable}.
\end{enumerate}
\end{thm}

\begin{table}
\centering
\begin{tabular}{|c|c|c|c|}
\hline
$\lambda$ & $\mathrm{hom}$ & $\mathrm{ext}^1$ & $\mathrm{ext}^2$ \\
\hline
\hline
$(1,0,0,0)$ & $1$ & $0$ & $1$ \\
\hline
$(1,1,0,0)$ & $1$ & $20$ & $2$ \\
\hline
$(2,0,0,0)$ & $1$ & $0$ & $191$ \\
\hline
$(2,1,0,0)$ & $1$ & $20$ & $401$ \\
\hline
$(2,1,1,0)$ & $1$ & $20$ & $191$ \\
\hline
$(2,2,0,0)$ & $1$ & $20$ & $590$ \\
\hline
$(3,0,0,0)$ & $1$ & $0$ & $5545$ \\
\hline
$(3,1,0,0)$ & $1$ & $20$ & $21419$ \\
\hline
$(3,1,1,0)$ & $1$ & $20$ & $10649$ \\
\hline
$(3,2,0,0)$ & $-$ & $-$ & $-$ \\
\hline
$(3,2,1,0)$ & $1$ & $-$ & $-$ \\
\hline
$(3,3,0,0)$ & $-$ & $-$ & $-$ \\
\hline
$(4,0,0,0)$ & $1$ & $0$ & $53065$ \\
\hline
$(4,1,0,0)$ & $1$ & $-$ & $-$ \\
\hline
$(4,1,1,0)$ & $1$ & $20$ & $141746$ \\
\hline
$(4,2,0,0)$ & $-$ & $-$ & $-$ \\
\hline
$(4,2,1,0)$ & $1$ & $-$ & $-$ \\
\hline
$(4,2,2,0)$ & $1$ & $20$ & $172910$ \\
\hline
$(4,3,0,0)$ & $-$ & $-$ & $-$ \\
\hline
$(4,3,1,0)$ & $-$ & $-$ & $-$ \\
\hline
$(4,4,0,0)$ & $-$ & $-$ & $-$ \\
\hline
\end{tabular}
\caption{Cohomological values for $\lambda_1 < 5$.}
\label[table]{MainTable}
\end{table}

\begin{rmk}
Note that the values appearing in \cref{MainTable} coincide with those computed in \cite[Theorem C]{FatOno} for $X_{BD}$. This appears not to be coincidental and suggests a deformation-theoretic correspondence between the same Schur functors on the two models of \hk manifolds, similar to the phenomenon described in \cite[Proposition 3.9]{OG26}. This correspondence is expected to hold for every Schur functor of the restricted quotient bundle, although no general proof is currently known. In particular, for $\bigwedge^2\mQ$, the agreement of the $\Ext^i$ dimensions was already observed in \cite[Remark 2.11]{Fat24}. The expected explanation is that the projectivization of the restricted quotient bundle on the Beauville--Donagi fourfold deforms to the projectivization of the corresponding bundle on the Debarre--Voisin fourfold.
\end{rmk}

\begin{rmk}
The computations on the Debarre--Voisin fourfold are more involved than those on the Beauville--Donagi fourfold, since the codimension of the former inside the ambient Grassmannian is $20$, compared to $4$ in the latter case. This results in a significantly longer Koszul complex, which is harder to handle computationally. In particular, indeterminacy (see \cref{def_indeterminate}) in the cohomology appears earlier, and we are able to compute exact values only for finitely many partitions. In \cref{MainTable}, we restrict to partitions $\lambda$ such that $\lambda_4 = 0$ and $\lambda_1 \geq \lambda_2 + \lambda_3$, since the graded algebra $\mathrm{Ext}^*(\Sigma_\lambda \mQ,\Sigma_\lambda \mQ)$ is invariant under tensoring with line bundles and under dualization. Moreover, we impose the bound $\lambda_1 < 5$, since for larger values the Koszul computations become intractable due to the proliferation of non-trivial terms. The symbol "$-$" in \cref{MainTable} indicates that the corresponding value cannot be computed with the available methods.
\end{rmk}

This approach to constructing modular sheaves with potentially non-trivial moduli spaces is not the only one. Indeed, O'Grady in \cite{OG26} exhibited examples of non-rigid modular sheaves on \hk manifolds of $\textrm{K3}^{[2]}$-type whose deformation spaces have dimension $2(a^2+1)$ for $a>1$. In future work, we aim to compare the techniques of \cite{OG26}, \cite{Fat24}, \cite{FatOno}, and the present paper.

Finally, we state a sufficient condition for the modular bundles constructed on $X_{DV}$ to be atomic:

\begin{thm}[\cref{atomicity}]
Let $X_{DV}\subseteq \mathrm{Gr}(6,10)$ be a very general Debarre--Voisin \hk manifold, and let $\lambda = (m,t,s,0)$ with $m\geq t+s$. Then $\Sigma_\lambda \mQ$ is atomic if and only if $\lambda=(m,0,0,0)$.
\end{thm}

\subsection*{Plan of the paper} The paper is divided as follows.
In \cref{sec:workTools} we recall some classical results about homogeneous varieties and homogeneous vector bundles. In \cref{sec:ChernClasses} we compute the Chern classes of any Schur functor of $\widetilde{\mQ}$ restricted to the Debarre--Voisin fourfold to prove that any such bundle is modular. 

In \cref{sec:cohom} we outline the general method applied to compute the dimensions of the $\Ext$-groups. Moreover, we state and prove a combinatorial result (\cref{LRDComputation}) that gives some sufficient conditions to get a moduli space with tangent space of positive dimension. 

In \cref{WorkedExample} and \cref{AboutSym}, we go through some worked out examples. In particular, in \cref{WorkedExample}, we compute explicitly the $\Ext$-groups of the first case with $\mathrm{ext}^1=20$, while in \cref{AboutSym}, we study the vector bundles $\Sym^m\mQ$ in more details. 

In the end, in \cref{sec:atom}, we study the atomicity of the Schur functors of $\mQ$. In particular, we characterize the partitions that provide atomic sheaves among the Schur functors of $\mQ$. 

\subsection*{Notations} Throughout this paper we work over $\mathbb{C}$ and $V_n$ refers to an $n$-dimensional complex vector space. We denote by $\Gr(k,n)$ the Grassmannian that parametrizes the $k$-dimensional subspaces $V_k$ of $V_n$. We denote by $\widetilde{\mU}$ and $\widetilde{\mQ}$, respectively, the tautological and the quotient bundle of $\Gr(k,n)$, fitting into the short exact sequence of vector bundles $$0\to\widetilde{\mU}\to\of_{\Gr(k,n)}\otimes V_n\to\widetilde{\mQ}\to0.$$ 
In what follows, we adopt the following notation: if $M$ is a complex variety and $\mathcal{F}$ is a vector bundle on $M$ we denote the zero locus of a global section $s$ of $\mathcal{F}$ by $X=\mZ(s)\subset M$. 

\subsubsection*{Acknowledgment} We would like to thank Valeria Bertini, Enrico Fatighenti, Francesco Meazzini, Kieran O'Grady, and Claudio Onorati for the insightful conversations and comments. We, moreover, thank Fabio Tanturri for sharing his knowledge about most of the tools needed for this work.\\ The authors acknowledge the partial support of the European Union - NextGenerationEU under the
National Recovery and Resilience Plan (PNRR) - Mission 4 Education and research - Component 2
From research to business - Investment 1.1 Notice Prin 2022 - DD N. 104 del 2/2/2022, from title
“Symplectic varieties: their interplay with Fano manifolds and derived categories”, proposal code
2022PEKYBJ – CUP J53D23003840006. The first author is partially supported by the Research Project PRIN 2020 - CuRVI, CUP J37G21000000001 and wishes to thank the MIUR Excellence Department
of Mathematics, University of Genoa, CUP D33C23001110001.\\
Both authors are members of the INDAM-GNSAGA group.

\section{Working tools}\label[section]{sec:workTools}

This section introduces the setting and all the tools we use to perform our computations. In particular, we will use the specific structure of the general Debarre--Voisin fourfold as a zero locus of a section inside a Grassmannian. Notice that the cohomological computations we are going to perform hold true for the general section of a homogeneous vector bundle on $\Gr(6,10)$, since we only need that the zero locus is of the expected dimension. 

However, in \cref{sec:ChernClasses}, we will need to consider a very general Debarre--Voisin fourfold, i.e., outside of the Noether--Lefschetz divisors, to ensure that the Picard rank is $1$.

\subsection{Homogeneous varieties}\label[section]{homogeneousV}
The general Debarre--Voisin fourfold can be seen as $$X_{DV} :=\mZ(s)\subset\Gr(6,10),$$ hence it corresponds to the zero locus of a general global section $$s\in H^0\left(\mathrm{Gr}(6,10),\bigwedge^3\widetilde{\mU}^{\vee}\right) $$ of a homogeneous vector bundle in a homogeneous variety. Let us recall some definitions, standard references are \cite{Ott95}, \cite{Sno89}, or \cite{Wey03}.

\begin{definition}
    Let $G$ be a connected, simply connected semi-simple complex Lie group and let $P$ be a parabolic subgroup of $G$. The quotient $G/P$ is a smooth projective variety and is called \textit{homogeneous variety}. A vector bundle $E$ on $G/P$ is called \textit{homogeneous} if the natural action of $G$ by multiplication on the base lifts to $E$. In other words, the following diagram commutes:
    \[
    \begin{tikzcd}
    G\times E \arrow[r] \arrow[d] & E \arrow[d] \\
    G\times G/P \arrow[r]         & G/P        
    \end{tikzcd}
    \]
\end{definition}
In our case, the Grassmannian $\Gr(6,10)$ can be realized as the quotient of $\mathrm{SL}(10)$ by the parabolic subgroup 
\[
P:= \left\lbrace \begin{pmatrix} A & * \\ 0 & B \end{pmatrix} \in \mathrm{SL}(10) ~ \mid ~ A \in \mathrm{SL}(6), ~B \in \mathrm{SL}(4)\right\rbrace.
\]

Recall that the irreducible representations of $\mathrm{SL}(N)$ are parametrized by partitions of length $N$. More precisely, let $\lambda\in\Z^{N}$, we call it partition of length $N$ if $\lambda_i\geq\lambda_{i+1}$ for $1\leq i\leq N$. Note that we do not require the entries to be non-negative.

The irreducible representation associated to $\lambda$ corresponds to the Schur functor $\Sigma_{\lambda}V_{10}$ applied to $V_{10}$. Its highest weight in the root system of $\mathfrak{sl}(10) :=Lie(\mathrm{SL}(10))$ is 
\[
[\lambda_1-\lambda_2,\lambda_2-\lambda_3,...,\lambda_{9}-\lambda_{10}].
\]
In particular, the partitions $\lambda$ and $\lambda+(a,...,a)$ give the same representation. Vice-versa, two Schur functors $\Sigma_{\lambda}V_{10}$ and $\Sigma_{\mu}V_{10}$ are isomorphic representations only if the differ by a partition of the form $(a,...,a)$. Thus, we can always assume that the entries of a partition are in $\mathbb{Z}_{\geq 0}$.

The representation associated with a partition $\lambda$ corresponds to the Schur functor $\Sigma_{\lambda}V_{N}$ applied to $V_N$, an $N$-dimensional vector space with the natural $\mathrm{SL}(N)$-action. In the same way, if $\lambda\in\mathbb{Z}^4$ and $\mu\in\mathbb{Z}^6$ are two partitions we denote
\[
    \Sigma_{\lambda|\mu}:=\Sigma_{\lambda}V_4\otimes\Sigma_{\mu}V_6
\]
the irreducible representation of $\mathrm{SL}(4)\times \mathrm{SL}(6)$ indexed by $(\lambda,\mu)$.\\
\begin{rmk}
The tautological bundle $\widetilde{\mU}$ on the Grassmannian $\Gr(6,10)$ is a homogeneous vector bundle of rank $6$. It corresponds to the standard representation of $\mathrm{SL}(6)$, while the quotient bundle $\widetilde{\mQ}$ is a homogeneous vector bundle of rank $4$. It corresponds to the standard representation of $\mathrm{SL}(4)$. 
\end{rmk}

On $\Gr(6,10)$, we still have a correspondence between (irreducible) homogeneous vector bundles and pairs of partitions. For example:
\begin{itemize}
\item $\lambda=(1,0,0,0)$ corresponds to the quotient bundle $\widetilde{\mQ}$;
\item $\mu=(1,0,0,0,0,0)$ corresponds to the tautological bundle $\widetilde{\mU}$;
\item $\lambda = (1,1,0,0)$ corresponds to $\bigwedge^2\widetilde{\mQ}$;
\item $\lambda = (k,0,0,0)$ corresponds to $\mathrm{Sym}^k\widetilde{\mQ}$;
\item for any pair of partitions $\lambda\in\mathbb{Z}^4$ and $\mu\in\mathbb{Z}^6$ we denote
\[\Sigma_{\lambda}\widetilde{\mQ}, \ \Sigma_{\mu}\widetilde{\mU} \]  
the irreducible homogeneous vector bundles corresponding to the Schur functors
\[
\Sigma_{\lambda|(0,0,0,0,0,0)}, \ \text{and} \ \Sigma_{(0,0,0,0)|\mu},\]
respectively.
\end{itemize}
With an abuse of notation, we write
\[
    \Sigma_{\lambda|\mu}:=\Sigma_{\lambda}\widetilde{\mQ}\otimes\Sigma_{\mu}\widetilde{\mU}.
\]
Some actions we can perform on the partitions reflect actions on the associated homogeneous bundle, and vice versa. For instance, if we call $\lambda+d:=(\lambda_1+d,...,\lambda_N+d)$ for any $d\in\mathbb{Z}$, then 
\[
\Sigma_{\lambda+d}V_N=\Sigma_{\lambda}V_N\otimes \left(\bigwedge^N V_N\right)^{\otimes d}\cong \Sigma_{\lambda}V_N,
\]

Moreover, the \textit{dual} Schur functor is given by
\[
\lambda^*:=(\lambda_1-\lambda_N,...,\lambda_1-\lambda_1).
\]

This implies that, with the above notation 

\[
\Sigma_{\lambda + d} \widetilde{\mQ} =\Sigma_{\lambda}\widetilde{\mQ} \otimes \mathcal{O}_{\Gr(6,10)}(d) \text{, ~~~~} \Sigma_{\mu+d}\widetilde{\mU} = \Sigma_{\lambda}\widetilde{\mU} \otimes \mathcal{O}_{\Gr(6,10)}(-d).
\]

In particular, since \[
\Sigma_{\lambda+d}\widetilde{\mQ}=\Sigma_{\lambda}\widetilde{\mQ}\otimes(\det\widetilde{\mQ})^{\otimes d}\ \text{and}\ \Sigma_{\mu+d}\widetilde{\mU}=\Schur_{\mu}\widetilde{\mU}\otimes(\det\widetilde{\mU})^{\otimes d}\] it follows that $\Sigma_{\lambda+d|\mu+d}=\Sigma_{\lambda|\mu}$.

\begin{rmk}\label[rmk]{CasiUtili}
    We can combine the above identities to get the following equality:
    \[
    \Schur_{(\lambda_1,\lambda_2,\lambda_3,\lambda_4)}\widetilde{\mQ}=\Schur_{(\lambda_1-\lambda_4,\lambda_2-\lambda_4,\lambda_3-\lambda_4,0)}\widetilde{\mQ}\otimes\of_{\mathrm{Gr}(6,10)}(\lambda_4).
    \]
    Since we are interested in the sheaf of endomorphisms of the Schur functor $\Schur_{(\lambda_1,\lambda_2,\lambda_3,\lambda_4)}\mQ$, we can restrict our analysis to partitions of the form $(m,t,s,0)$.
    
    Moreover, since 
    \[
\End(\Schur_{\lambda}\widetilde{\mQ})\cong\Schur_{\lambda}\widetilde{\mQ}\otimes\Schur_{\lambda}\widetilde{\mQ}^{\vee},
    \] 
    it is equivalent to consider the partition $(m,t,s,0)$ or the partition associated to the dual vector bundle $(m,m-s,m-t,0)$. For this reason, we restrict ourselves to partitions $(m,t,s,0)$ with $m\geq t+s$. Indeed, if $m<s+t$, then $m>2m-(s+t)=(m-s)+(m-t)$, which means that the inequality holds for $(m,m-s,m-t,0)$ in place of $(m,t,s,0)$.
\end{rmk}

\begin{notz}\label[notz]{notz:endBundle}
    For any non--increasing sequence $\lambda\in\mathbb{Z}^4$ we consider the Schur functor $\Sigma_{\lambda}\widetilde{\mQ}$ and we denote by
\[
\widetilde{\mathcal{E}}_\lambda:=\mathcal{E}nd(\Sigma_{\lambda}\widetilde{\mQ}) =\Sigma_{\lambda}\widetilde{\mQ}\otimes\Sigma_{\lambda}\widetilde{\mQ}^{\vee}
\]
the associated endomorphism bundles. When we restrict to the Debarre--Voisin fourfold $X$, we write 
\[
    \mQ := \widetilde{\mQ}|_X, \   \ \Sigma_{\lambda}\mQ:=(\Sigma_{\lambda}\widetilde{\mQ})|_X\ \text{and}\ \mathcal{E}_{\lambda}:=(\widetilde{\mathcal{E}}_{\lambda})|_X.
\]
\end{notz}

\begin{rmk}
A fundamental tool in this setting is the Borel--Weil--Bott theorem (see, for instance, \cite[Corollary 4.1.9]{Wey03}), which allows one to compute the cohomology of any irreducible homogeneous vector bundle by considering only the associated partitions. This will be particularly useful, as it enables these computations to be carried out algorithmically. In the following lemma we write the Borel--Weil--Bott theorem for $\Gr(6,10)$.
\end{rmk}

\begin{lemma}\label{lm:BB} \cite[Corollary 4.1.9]{Wey03}
Let $G=\Gr(6,10)$, let $\lambda$ be a partition of length $4$, and let $\mu$ be a partition of length $6$. Set $\gamma=(\lambda,\mu)$. Then exactly one of the following mutually exclusive cases occurs:
\begin{itemize}
    \item There exists a nontrivial permutation $\sigma\in S_{10}$ such that
    \[
    \sigma(\gamma+(9,\dots,1,0))-(9,\dots,1,0)=\gamma.
    \]
    In this case, $H^i(G,\Sigma_\gamma)=0$ for all $i\geq 0$.
    
    \item There exists a permutation $\sigma\in S_{10}$ such that
    \[
    \sigma(\gamma+(9,\dots,1,0))-(9,\dots,1,0)=\tau
    \]
    is a partition. Then
    \[
    H^i(G,\Sigma_\gamma)=
    \begin{cases}
        0 & \text{for } i\neq \ell(\sigma),\\
        \Sigma_{\tau} V_{10} & \text{for } i=\ell(\sigma),
    \end{cases}
    \]
    where $\ell(\sigma)$ denotes the length of $\sigma$ (i.e.\ the number of inversions).
\end{itemize}
\end{lemma}

\subsection{Combinatorial rules}
As we have seen, partitions and associated Schur functors play an important role in this setting. In particular, to compute the tensor product of Schur functors of vector bundles, we need to understand how to multiply Schur functors. This is given by the Littlewood--Richardson's rule, see e.g. \cite[Theorem 2.3.4]{Wey03} or \cite[Theorem 4.9.4]{Sag}.
\begin{alg}[Littlewood--Richardson's rule]\label[alg]{LRD} 
    Let $\lambda$, $\mu$ be two partitions of the same length $n$. The product of the associated Schur functors decomposes as 
    \[
    \Schur_\lambda \otimes \Schur_\mu = \bigoplus_\nu N_{\lambda,\mu}^\nu \Schur_\nu.
    \]
    The sum runs through all the partitions $\nu$ of $|\lambda|+|\mu|$ of length $n$ and $N_{\lambda,\mu}^\nu$ are the \textit{Littlewood--Richardson's numbers}, which can be computed in the following way. \\
    We start with the Young diagram associated with the partition $\lambda$ (in gray), and we add $|\mu|$ blocks in a way such that the result is again a Young diagram corresponding to a partition $\nu$.
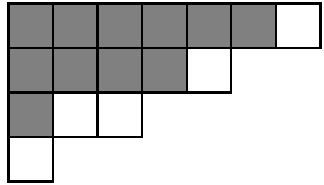
    \[
    \ytableausetup{centertableaux}
    \ytableaushort
    {\none\none\none \none\none ,\none\none\none, \none}
    * {7,5,3,1}
    * [*(gray)]{6,4,1,0}
    \]
    Now we fill the blocks we added (in white) with the \textit{content} of $\mu$, i.e., with $\mu_1$ "ones", $\mu_2$ "twos"...  following three rules:
    \begin{itemize}
        \item[(a)] for each row, the row-word has to be weakly increasing (left to right);
        \item[(b)] for each column, the column-word has to be strictly increasing (top to down);
        \item[(c)] starting from top-right, if we list the entries row by row, from right to left, at each step the number of "ones" has to be greater or equal to the number of "twos", the number of "twos" has to be greater or equal to the number of "threes" and so on. 
    \end{itemize}
    The number of admissible fillings of the shape $\nu$ corresponds to $N_{\lambda,\mu}^\nu$. 
\end{alg}

\begin{ex}
    Let $\lambda = (2,1,0)$ and $\mu = (2,2,0)$ then Littlewood--Richardson's rule returns
    \[
    \ytableausetup{boxsize=1.25em}
\ytableausetup{aligntableaux=top}
\ydiagram[*(gray)]{2,1} ~~~~
\otimes ~~~~\ydiagram{2,2}~~~~
=~~~~\ytableaushort{\none\none 1 1, \none 2 2}
    * {4,3}
    * [*(gray)]{2,1}
\oplus  ~~~~
\ytableaushort{\none\none 1 1, \none 2, 2}
    * {4,2,1}
    * [*(gray)]{2,1}
\oplus  ~~~~
\ytableaushort{\none\none 1 ,\none 1 2 , 2}
    * {3,3,1}
    * [*(gray)]{2,1}
    \oplus ~~~~
    \ytableaushort{\none\none 1 ,\none 1, 2 2}
    * {3,2,2}
    * [*(gray)]{2,1}
\]
\end{ex}

\begin{rmk}
A special case of Littlewood--Richardson's rule occurs when $\mu=(m,0,...,0)$, i.e. when we tensor with a symmetric power of a vector bundle. This particular case goes under the name of \textit{Pieri's rule}. In this case, since we have to label the white blocks with $m$ ones, Littlewood--Richardson's numbers are either zeros or ones. The admissible partitions $\nu$ are the ones obtained by $\lambda$ by adding $m$ blocks but no more than one for each column.

\begin{ex}
    Let $\lambda = (2,1,0)$ and $m = 3$ then Pieri's rule gives
    \[
    \ytableausetup{boxsize=1.25em}
\ytableausetup{aligntableaux=top}
\ydiagram[*(gray)]{2,1} ~~~~
\otimes ~~~~\ydiagram{3}~~~~
=~~~~\ytableaushort{\none\none , \none }
    * {5,1}
    * [*(gray)]{2,1}
\oplus  ~~~~
\ytableaushort{\none\none,\none }
    * {4,2}
    * [*(gray)]{2,1}
    \oplus ~~~~
    \ytableaushort{\none\none,\none }
    * {4,1,1}
    * [*(gray)]{2,1}
    \oplus ~~~~
    \ytableaushort{\none\none ,\none}
    * {3,2,1}
    * [*(gray)]{2,1}
\]
\end{ex}
\end{rmk}

\subsection{Koszul complex.} \label[section]{CohomologyComputations}
The last tool we want to introduce is the Koszul complex. Consider a smooth projective variety $M$ and a globally generated vector bundle $\widetilde{\mF}$ of rank $r$ on $M$. Let $X: =\mZ(s)$, where $s\in H^0\left(M,\widetilde{\mF}\right)$, so that the normal bundle of $X$ in $M$ corresponds to $\mF:=\mathcal{\widetilde{\mF}}_{|X}$. If we denote by $\iota:X\hookrightarrow M$ the embedding, under these hypothesis the Koszul complex gives the following resolution of $\of_X$ (or, to be more precise, of the $\of_M$-module $\iota_{\ast}\of_X$):
\[
0\to \mathrm{det}\left(\widetilde{\mF}^{\vee}\right)\to \bigwedge^{r-1}\widetilde{\mF}^{\vee}\to...\to\bigwedge^2\widetilde{\mF}^{\vee}\to\widetilde{\mF}^{\vee}\to\of_M\to \iota_*\of_X\to 0.
\]
We specialize to the case $M= \Gr(6,10)$ and $\Tilde{F}=\bigwedge^3\widetilde{\mU}^{\vee}$. The Koszul complex of $\iota: X\hookrightarrow \Gr(6,10)$ is:
\begin{equation}\label{Koszul}
0\to\det\left(\bigwedge^3\widetilde{\mU}\right)\to\bigwedge^{19}\bigwedge^3\widetilde{\mU}\to...\to\bigwedge^2\bigwedge^3\widetilde{\mU}\to\bigwedge^3\widetilde{\mU}\to\of_{\Gr(6,10)}\to \iota_*\of_X\to0.
\end{equation}
\begin{rmk}\label[rmk]{rmk:dWtoSchur}
By \cite[Ex. 6.17]{FH91} we can write each term $\bigwedge^p\bigwedge^3\Tilde\mU$ as direct sums of $\Schur_{\mu}\Tilde\mU$. This can be done via computer algebra software such as \red{Macaulay2}. Moreover, notice that if $\widetilde{\mU}$ is the tautological bundle of $\Gr(6,10)$ and $k\in\mathbb{Z}_{\geq 1}$ then $$\bigwedge^{10+k}\bigwedge^3\widetilde{\mU}\cong \left(\bigwedge^{10-k}\bigwedge^3\widetilde{\mU}\right)\otimes\of(-k).$$
To see this, one has to apply the formula
\[
    \bigwedge^{\mathrm{rk}F-k}F\cong \bigwedge^kF^{\vee}\otimes\det F
\]
which holds for every vector bundle $F$ (see e.g \cite[Ex. 5.16(b)]{Har77}) twice, namely to $F=\widetilde\mU$ and to $F=\bigwedge^3\widetilde\mU$.
\end{rmk}

Using \cref{rmk:dWtoSchur}, we decompose all the bundles appearing in \eqref{Koszul} into sums of irreducible factors. We stored all the information in \cref{tab:your_label}.

\begin{table} 
    \centering
    \begin{adjustbox}{max width=\textwidth}
    \begin{tabular}{|c|c|c|c|c|c|c|c|c|c|c|}
        \hline
         $p=0$ & $p=1$ & $p=2$ & $p=3$ & $p=4$ & $p=5$ & $p=6$ & $p=7$ & $p=8$ & $p=9$ & $p=10$ \\ 
        \hline
         (0,0,0,0,0,0) & (1,1,1,0,0,0) & (2,2,1,1,0,0) & (3,3,1,1,1,0) & (4,4,1,1,1,1) & (5,4,2,2,1,1) & (6,4,3,2,2,1) & (7,4,4,2,2,2) & (8,4,4,3,3,2) & (9,4,4,4,3,3) & (10,4,4,4,4,4)\\ 
        \hline
         - & - & (1,1,1,1,1,1) & (3,2,2,2,0,0) & (4,3,2,2,1,0) & (5,3,3,2,2,0) & (6,3,3,3,3,0) & (7,4,3,3,3,1) & (7,5,5,3,2,2) & (8,5,5,4,3,2) & (9,5,5,5,3,3) \\ 
        \hline
         - & - & - & (2,2,2,1,1,1) & (3,3,3,3,0,0) & (4,4,3,3,1,0) & (5,5,3,3,1,1) & (6,5,4,3,2,1) & (7,5,4,4,3,1) & (8,5,5,3,3,3) & (9,5,5,4,4,3) \\ 
        \hline
         - & - & - & - & (3,3,2,2,1,1) & (4,4,2,2,2,1) & (5,5,2,2,2,2) & (6,5,3,3,2,2) & (7,5,4,3,3,2) & (8,5,4,4,4,2) & (8,6,6,4,4,2) \\ 
        \hline
         - & - & - & - & (2,2,2,2,2,2) & (4,3,3,3,1,1) & (5,4,4,3,2,0) & (6,4,4,4,3,0) & (7,4,4,4,4,1) & (7,6,6,3,3,2) & (8,6,6,4,3,3) \\ 
        \hline
         - & - & - & - & - & (3,3,3,2,2,2) & (5,4,3,3,2,1) & (6,4,4,3,3,1) & (6,6,5,3,3,1) & (7,6,5,4,4,1) & (8,6,5,5,4,2) \\ 
        \hline
         - & - & - & - & - & - & (4,4,4,4,1,1) & (5,5,5,3,3,0) & (6,6,4,4,2,2) & (7,6,5,4,3,2) & (8,6,5,4,4,3) \\ 
        \hline
         - & - & - & - & - & - & (4,4,3,3,2,2) & (5,5,4,4,2,1) & (6,6,3,3,3,3) & (7,6,4,4,3,3) & (8,5,5,5,5,2) \\ 
        \hline
         - & - & - & - & - & - & (3,3,3,3,3,3) & (5,5,3,3,3,2) & (6,5,5,4,4,0) & (7,5,5,5,4,1) & (7,7,7,3,3,3) \\ 
        \hline
         - & - & - & - & - & - & - & (5,4,4,4,2,2) & (6,5,5,4,3,1) & (7,5,5,4,4,2) & (7,7,6,4,4,2) \\ 
        \hline
         - & - & - & - & - & - & - & (4,4,4,3,3,3) & (6,5,4,4,3,2) & (6,6,6,4,4,1) & (7,7,5,5,5,1) \\ 
        \hline
         - & - & - & - & - & - & - & - & (5,5,5,5,2,2) & (6,6,5,5,5,0) & (7,7,5,5,3,3) \\ 
        \hline
         - & - & - & - & - & - & - & - & (5,5,4,4,3,3) & (6,6,5,5,3,2) & (7,7,4,4,4,4) \\ 
        \hline
         - & - & - & - & - & - & - & - & (4,4,4,4,4,4) & (6,6,4,4,4,3) & (7,6,6,5,5,1) \\ 
        \hline
         - & - & - & - & - & - & - & - & - & (6,5,5,5,3,3) & (7,6,6,5,4,2) \\ 
        \hline
         - & - & - & - & - & - & - & - & - & (5,5,5,4,4,4) & (7,6,5,5,4,3) \\ 
        \hline
         - & - & - & - & - & - & - & - & - & - & (6,6,6,6,6,0) \\ 
        \hline
         - & - & - & - & - & - & - & - & - & - & (6,6,6,6,3,3) \\ 
        \hline
         - & - & - & - & - & - & - & - & - & - & (6,6,5,5,4,4) \\ 
        \hline
         - & - & - & - & - & - & - & - & - & - & (5,5,5,5,5,5) \\ 
        \hline
    \end{tabular}
    \end{adjustbox}
    \caption{Irreducible factors in the Koszul complex of $X_{DV}$.}
    \label[table]{tab:your_label}
\end{table}

\section{Chern classes}\label[section]{sec:ChernClasses}
In this section we explicitly compute the Chern classes of any Schur functor $\Sigma_{m,t,s,0} \mQ$ and hence its discriminant $\Delta(\Sigma_{m,t,s,0} \mQ)$. This explicit computation will prove item (1) of \cref{MainTheorem}. 

We start by recalling some well-known facts for \hk manifolds of K3$^{[2]}$-type (see e.g. \cite[Lemma 1.6]{FatOno}).
\begin{lemma}\label[lemma]{toddIHS}
    Let $X$ be a \hk manifold of $\mathrm{K3}^{[2]}$-type. Denote by $\mathsf{p}\in H^8(X,\mathbb{C})$ the class of a point. Then the following relations hold:
    \begin{enumerate}
        \item $\int_X \mathsf{c}_2(X)^2 = 828$;
        \item $\int_X \mathsf{c}_4(X) = 324$;
        \item $\mathrm{td}_X= 1 + \frac{1}{12} \mathsf{c}_2(X) +3\mathsf{p}$;
        \item $\sqrt{\mathrm{td}_X}= 1 + \frac{1}{24} \mathsf{c}_2(X) +\frac{25}{32}\mathsf{p}$.
    \end{enumerate}
\end{lemma}

The following lemma tells us the relations between the Chern characters of $\mQ$ coming from the intersection ring of a very general Debarre--Voisin fourfold $X$. Note that, since the Picard number of $X$ is $1$, we have an isometric isomorphism between the Picard lattice of the Grassmannian $\mathrm{Gr}(6,10)$ and $H^{1,1}(X,\Z)$. These results will allow us to express the Chern characters of the Schur functor $\Sigma_{m,t,s,0} \mQ$ in terms of the Chern characters of $\mQ$.
\begin{lemma}\label[lemma]{ChernRelations}
    Let $X\subseteq \mathrm{Gr}(6,10)$ be a very general Debarre--Voisin \hk manifold, then:
    \begin{enumerate}
        \item $\Delta(\mQ)=\ch_1(\mQ)^2 -8\,\ch_2(\mQ)=\mathsf{c}_2(X)$.
        \item The group $\mathrm{H}^{3,3}(X,\mathbb{Q})$ has rank $1$ and the following relations hold:
        \begin{enumerate}
            \item $\ch_1(\mQ)^3=-264\,\ch_3(\mQ)$;
            \item $\ch_1(\mQ)\cdot \ch_2(\mQ)=-18\,\ch_3(\mQ)$.
        \end{enumerate}
        \item The group $\mathrm{H}^{8}(X,\mathbb{Q})$ has rank $1$ and the following relations hold:
        \begin{enumerate}
            \item $\ch_1(\mQ)^4=-5808\,\ch_4(\mQ)$;
            \item $\ch_1(\mQ)^2\cdot \ch_2(\mQ)=-396\,\ch_4(\mQ)$;
            \item $\ch_1(\mQ)\cdot \ch_3(\mQ)=22\,\ch_4(\mQ)$;
            \item $\ch_2(\mQ)^2=-60\, \ch_4(\mQ)$;
            \item $\int_X\ch_4(\mQ)=-\frac{1}{4}$.
            \end{enumerate}
        \end{enumerate}
    In particular, $\mQ$ is modular.
\end{lemma}
\begin{proof}
    These relations follow from the \textit{normal bundle sequence}
    \[
    0\longrightarrow \mathcal{T}_X \longrightarrow \mU^\vee\otimes \mQ \longrightarrow \bigwedge^3 \mU^\vee \longrightarrow 0.
    \]
    Following the notations in \cite{DV10}, we denote $\mathsf{c}_i(\mU^\vee)$ by $\mathsf{c}_i$. Using the tautological sequence of the Grassmannian $\mathrm{Gr}(6,10)$ we get
    \begin{align*}
        \ch_1(\mQ) &= \mathsf{c}_1,\\
        \ch_2(\mQ) &= \mathsf{c}_2- \frac{1}{2}\mathsf{c}_1^2,\\
        \ch_3(\mQ) &= \frac{1}{2}\left(\mathsf{c}_3-\mathsf{c}_1\mathsf{c}_2+\frac{1}{3}\mathsf{c}_1^3\right),\\
        \ch_4(\mQ) &= \frac{1}{6}\left(\mathsf{c}_4-\mathsf{c}_1\mathsf{c}_3+\mathsf{c}_1^2\mathsf{c}_2-\frac{1}{2}\mathsf{c}_2^2-\frac{1}{4}\mathsf{c}_1^4\right),
    \end{align*}
    and
    \[
    \ch (\mU^\vee\otimes \mQ) =24 + 10\ch_1(\mQ)+ \left(\ch_1(\mQ)^2 + 2\ch_2(\mQ)\right) + 10\ch_3(\mQ) + \left( 2\ch_4(\mQ)+2\ch_1(\mQ)\ch_3(\mQ) -\ch_2(\mQ)^2\right).
    \]
    By applying the \textit{splitting principle} to a rank $6$ vector bundle, we obtain
    \begin{align*}
    \ch \left(\bigwedge^3 \mU^\vee \right) = 20 &+ 10\ch_1(\mQ)+ \left(2\ch_1(\mQ)^2 -6\ch_2(\mQ)\right) + \left( \frac{1}{6}\ch_1(\mQ)^3-3\ch_1(\mQ)\ch_2(\mQ)\right) + \\ 
    &+\left( 6\ch_4(\mQ)+\ch_1(\mQ)\ch_3(\mQ) +\ch_2(\mQ)^2-\frac{1}{2}\ch_1(\mQ)^2\ch_2(\mQ)\right).
    \end{align*}
    On the other hand, for a \hk manifold of type $\textrm{K3}^{[2]}$ we have $ \mathsf{c}_2(X)^2 = 828\mathsf{p}$ and $ \mathsf{c}_4(X) = 324\mathsf{p}$. We deduce 
    \[
    \ch(\mathcal{T}_X) = 4 + 0- \mathsf{c}_2(X)+0 + 15\mathsf{p}.
    \]
    Hence, we have the following equalities
    \begin{align*}
        \mathsf{c}_2(X) &= \ch_1(\mQ)^2 - 8\ch_2(X),\\
        10\ch_3(\mQ) &= \frac{1}{6}\ch_1(\mQ)^3-3\ch_1(\mQ)\ch_2(\mQ),\\
        15\mathsf{p} &=  -4\ch_4(\mQ)+\ch_1(\mQ)\ch_3(\mQ) -2\ch_2(\mQ)^2+\frac{1}{2}\ch_1(\mQ)^2\ch_2(\mQ).
    \end{align*}
    The first one yields (1). From the intersection numbers
    \[
    \mathsf{c}_1\mathsf{c}_3=330\mathsf{p}~~~~~ \mathsf{c}_4=105\mathsf{p}~~~~~ \mathsf{c}_1^2\mathsf{c}_2=825\mathsf{p}~~~~~
    \mathsf{c}_2^2=477\mathsf{p}~~~~~
    \mathsf{c}_1^4=1452\mathsf{p},
    \]
    (see the proof of \cite[Eq. 10]{DV10}) we get 
    \[
    \mathsf{ch}_1(Q)\mathsf{ch}_3(Q)=-\frac{11}{2}\mathsf{p}~~~ \mathsf{ch}_4(Q)=-\frac{1}{4}\mathsf{p}~~~ \mathsf{ch}_1(Q)^2\mathsf{ch}_2(Q)=99\mathsf{p}~~~
    \mathsf{ch}_2(Q)^2=15\mathsf{p}~~~
    \mathsf{ch}_1(Q)^4=1452\mathsf{p}.
    \]
    Part (3) of the assertion follows from this.\\
    Since $X$ is a very general projective \hk manifold of $\textrm{K3}^{[2]}$--type, its Picard number is $1$ and the group $H^{1,1}(X,\mathbb{Q})$ is generated by $\ch_1(\mQ)$. Dually, the group $H^{3,3}(X,\mathbb{Q})$ is generated by $\ch_1(\mQ)^\vee$, i.e., the class satisfying $\int_X\ch_1(\mQ)^\vee\cdot\alpha=q_X(\ch_1(\mQ),\alpha)$ for every $\alpha \in H^2(X,\mathbb{Q})$, where $q_X$ is the Beauville--Bogomolov--Fujiki form. By the Beauville--Bogomolov--Fujiki relation applied to $\ch_1(\mQ)$, we have $q_X(\ch_1(\mQ))=22$ and
    \[
    \ch_1(\mQ)^\vee \ch_1(\mQ) = 22\mathsf{p}.
    \]
    It follows that 
    \[
    \mathsf{ch}_1(\mQ)^3=66 \ch_1(\mQ)^\vee~~~  \mathsf{ch}_1(\mQ)\mathsf{ch}_2(\mQ)=\frac{9}{2} \ch_1(\mQ)^\vee~~~
    \mathsf{ch}_3(\mQ)=-\frac{1}{4} \ch_1(\mQ)^\vee,
    \]
    from which we get (2).
\end{proof}

We are now ready to compute the Chern characters of the general Schur functor $\Sigma_{m,t,s,0} \mQ$ and so its discriminant. To state the next results, we need the following notations: 
\begin{align*}
    r(m,t,s) &:= \frac{(m+3)(t+2)(s+1)(m-t+1)(m-s+2)(t-s+1)}{12},\\
    \ell (m,t,s) &:= \frac{m+t+s}{4},\\
    \delta(m,t,s)&:= \frac{3m^2-2mt-2ms+3t^2+3s^2-2ts+12m+4t-4s}{60},\\
    \tau(m,t,s)&:= 15\delta(m,t,s)-44\ell(m,t,s)^2,\\
    \alpha_3(t,s)&:= -60t-60s+30,\\
    \alpha_2(t,s)&:= -109t^2-241ts-109s^2+103t+80s-21,\\
    \alpha_1(t,s)&:= -60t^3-241t^2s-241ts^2-60s^3+65t^2+78ts+4s^2-t+8s+6,\\
    \alpha_0(t,s)&:= -10t^4-60t^3s-109t^2s^2-60ts^3-10s^4+10t^3+19t^2s-19ts^2-10s^3+\!3t^2-\!13ts+\!3s^2+\!14t-\!14s,\\
    \xi(m,t,s) &:= \frac{-10m^4+\alpha_3(t,s)m^3+\alpha_2(t,s)m^2+\alpha_1(t,s)m+\alpha_0(t,s)}{20}.
\end{align*}

\begin{rmk}\label[remark]{rmk_ChClasses}
    By \cref{ChernRelations}, the geometric relations in $H^*(X,\mathbb{Q})$ start appearing from degree $3$ on. This means that the computations the Chern characters $\ch_i$ for $i=0,1,2$ can be performed in $H^*(\mathrm{Gr}(6,10),\mathbb{Q})$, i.e., before the restriction. 

    Closed formulae for the Chern characters (in degrees up to $3$) of Schur functors appeared recently in \cite[Appendix B]{DFO25a}. In particular, we can use these explicit formulae to compute $\ch_i(\Sigma_\lambda \mQ)$ for $i=0,1,2,3$ as a polynomial in the Chern characters of $\mQ$. As mentioned before, for $i=3$, we also need to impose the relations stated in \cref{ChernRelations}.

    Clearly, in low degrees, closed formulae were already known: the rank of a Schur functor is computed through Weyl's dimension formula for representations of Lie groups of type $A_3$; the first Chern class has been established in \cite{Rub13}.
\end{rmk}

\begin{proposition}\label[proposition]{CompChClasses}
    Let $X\subseteq \mathrm{Gr}(6,10)$ be a very general Debarre--Voisin \hk manifold, and let $\lambda = (m,t,s,0)$ with $m\geq t+s$, then:
    \begin{enumerate}
        \item $\ch_0(\Sigma_\lambda \mQ)= r(m,t,s)$, which is the rank of $\Sigma_{m,t,s,0} \mQ$.;
        \item $\ch_1(\Sigma_\lambda \mQ)=\ell (m,t,s) r(m,t,s) \ch_1(\mQ)$;
        \item $\ch_2(\Sigma_\lambda \mQ)= \delta (m,t,s) r(m,t,s)\ch_2(\mQ) + \frac{1}{2}\left( \ell(m,t,s)^2-\frac{\delta(m,t,s)}{4}\right)r(m,t,s)\ch_1(\mQ)^2$;
        \item $\ch_3(\Sigma_\lambda \mQ)= \tau(m,t,s)\ell(m,t,s)r(m,t,s)\ch_3(\mQ)$;
        \item $\ch_4(\Sigma_\lambda \mQ)= \xi(m,t,s)r(m,t,s)\ch_4(\mQ)$.
        
    \end{enumerate}
\end{proposition}
\begin{proof} 
    By \cref{rmk_ChClasses}, we only need to compute $\ch_4(\Sigma_\lambda \mQ)$. The polynomial appearing in the statement is obtained by interpolating several values of $m$, $t$, and $s$ and by checking it by induction on the entries of the partition. 
    \begin{itemize}
        \item \underline{Case $t=0$} This corresponds to the case of the symmetric power and follows from \cref{ChernSym}.\\
        \item \underline{Case $t>0$, $s=0$} Here we use Pieri's rule on
        \[
        \mathrm{Sym}^m\mQ \otimes \mathrm{Sym}^t \mQ= \bigoplus_{i=0}^t \Schur_{m+t-i,i}\mQ
        \]
        in order to obtain the Chern character of $\Schur_{m,t}\mQ$ as
        \[
        \ch\left(\Schur_{m,t}\mQ\right) = \ch\left(\mathrm{Sym}^m\mQ\right) \ch\left(\mathrm{Sym}^t\mQ\right) - \sum_{i=0}^{t-1}\ch\left(\Schur_{m+t-i,i}\mQ\right).
        \]
        If we isolate the degree $d=2$, $3$ and $4$ we can use the inductive hypothesis to get $\ch_d\left(\Schur_{m,t}\mQ\right)$.
        \item \underline{Case $s>0$} Here we use Pieri's rule on $\Schur_{m,t}\mQ \otimes \mathrm{Sym}^s \mQ$, which gives 
        \begin{align*}
            \Schur_{m,t}\mQ \otimes \mathrm{Sym}^s \mQ = \Schur_{m+s,t}\mQ \oplus \Schur_{m+s-1,t+1}\mQ \oplus \Schur_{m+s-2,t+2}\mQ \oplus \dots &\oplus \Schur_{m,t+s}\mQ \oplus \\
            \oplus \Schur_{m+s-1,t,1}\mQ \oplus \Schur_{m+s-2,t+1,1}\mQ \oplus \dots &\oplus \Schur_{m,t+s-1,1}\mQ \oplus \\
             \oplus \Schur_{m+s-2,t,2}\mQ \oplus \dots &\oplus \Schur_{m,t+s-2,2}\mQ \oplus\\
             &\vdots\\
             \oplus\Schur_{m+1,t,s-1}\mQ &\oplus\Schur_{m,t+1,s-1}\mQ \\
             &\oplus\Schur_{m,t,s}\mQ.
        \end{align*}
        To perform induction, one applies the Chern character on both sides and extracts the term $\ch\left(\Schur_{m,t,s}\mQ\right)$. As before, for $d=2$, $3$ and $4$ we can use the inductive hypothesis to get $\ch_d\left(\Schur_{m,t,s}\mQ\right)$.
    \end{itemize}
    Then one has to check several equalities between long polynomials in several variables, which is straightforward but quite cumbersome. For this reason, we made the computations using computer algebra software, such as \red{Macaulay2} and \red{Wolfram Alpha}. 
\end{proof}

\begin{rmk}
    The polynomial appearing in the fourth Chern character is quite complicated. Unlike the previous ones, it is not possible to express it as a polynomial in $\delta$ and $\ell$. More precisely, it is not a linear combination of $\ell^4$, $\ell^2\delta$ and $\delta^2$. This phenomenon seems to be an instance of the difficulty of finding a closed formula in degree four; see the discussion in \cite[Appendix A]{DFO25a}.
\end{rmk}

We obtain two key results from \cref{CompChClasses}.
\begin{corollary}\label[corollary]{SchurIsModular}
Let $X\subseteq \mathrm{Gr}(6,10)$ be a very general Debarre--Voisin \hk manifold, and let $\lambda = (m,t,s,0)$ with $m\geq t+s$, then  \[\Delta(\Sigma_\lambda \mQ)= \frac{\delta (m,t,s)}{4} r(m,t,s)^2\mathsf{c}_2(X).\] 
In particular, $\Sigma_\lambda \mQ $ is modular.
\end{corollary}
\begin{proof}
    By Items (2) and (3) of \cref{CompChClasses} we get 
    \begin{align*}
    \Delta(\Sigma_\lambda \mQ) :=& \,\ch_1(\Sigma_\lambda \mQ)^2 - 2 \ch_0(\Sigma_\lambda \mQ)\ch_2(\Sigma_\lambda \mQ) = \\
    =& \,\ell^2 r^2 \ch_1(\mQ)^2 - 2r^2\delta\ch_2(\mQ) -\ell^2 r^2 \ch_1(\mQ)^2 +\frac{r^2\delta}{4}\ch_1(\mQ)^2 =\\
    =&\,\frac{r^2\delta}{4}\left(\ch_1(\mQ)^2 -8 \ch_2(\mQ)\right).
    \end{align*}
     
\end{proof}
\begin{rmk}
For a vector bundle $\mathcal{E}$, we have 
\[
\ch(\mathcal{E}nd(\mathcal{E})) = (\mathrm{rk}(\mathcal{E})^2, 0 , -\Delta(\mathcal{E}),0,\ch_2(\mathcal{E})^2-2\ch_1(\mathcal{E})\ch_3(\mathcal{E})+2\mathrm{rk}(\mathcal{E})\ch_4(\mathcal{E})).
\]
For this reason, we denote the top degree part by $\Xi(\mathcal{E}):= \ch_2(\mathcal{E})^2-2\ch_1(\mathcal{E})\ch_3(\mathcal{E})+2\mathrm{rk}(\mathcal{E})\ch_4(\mathcal{E})$, so that 
\[
\ch(\mathcal{E}nd(\mathcal{E})) = (\mathrm{rk}(\mathcal{E})^2, 0 , -\Delta(\mathcal{E}),0,\Xi(\mathcal{E})).
\]
\end{rmk}
\begin{corollary}\label[corollary]{SchurIsModular2}
Under the hypothesis of \cref{CompChClasses}, we have:
    \begin{enumerate}
        \item $\Xi(\Sigma_\lambda \mQ)= -\frac{1}{3312}\left(-m^4+\frac{\alpha_3}{10}m^3+\frac{\alpha_2}{10}m^2+\frac{\alpha_1}{10}m+\frac{\alpha_0}{10}+484\ell^4-330\ell^2\delta-\frac{207}{4}\delta^2\right)r^2\mathsf{c}_2(X)^2$;\\
        \item $\int_X\Xi(\Sigma_\lambda \mQ)=-\frac{1}{4}\left(-m^4+\frac{\alpha_3}{10}m^3+\frac{\alpha_2}{10}m^2+\frac{\alpha_1}{10}m+\frac{\alpha_0}{10}+484\ell^4-330\ell^2\delta-\frac{207}{4}\delta^2\right)r^2$;\\
        \item $\chi(\Sigma_\lambda \mQ,\Sigma_\lambda \mQ)= 3\left(1+\frac{-276\delta-1936\ell^4+1320\delta\ell^2+207\delta^2-8\xi}{48}\right)r^2$.
    \end{enumerate}
\end{corollary}
\begin{proof}
These formulae are obtained by applying \cref{CompChClasses} in the definition of $\Xi$ and in the Riemann–-Roch formula. Since the computations are quite cumbersome, we used \red{Wolfram Alpha}.
\end{proof}
\begin{rmk}\label[rmk]{analogyForChern}
    We note that the Chern characters of degree $0$, $1$ and $2$ are given by the same polynomials as in \cite[Theorem 4.1]{FatOno}, while they differ in degree $3$ and $4$. More precisely, the polynomials $r$, $\ell$ and $\delta$ are the same, while $\tau$ and $\xi$ are different. This is a consequence of the fact that the relations of the intersection ring of $X$ appear from degree $3$ onward. Nevertheless, by expanding the polynomials which give the Euler characteristic here and in the case of \cite{FatOno}, one can see that they are in fact the same. This will be useful in \cref{proof7.1}.
\end{rmk}

\section{Cohomological computations}\label[section]{sec:cohom}
In this section, we show the computation of $\mathrm{Ext}\left(\Sigma_\lambda \mQ,\Sigma_\lambda \mQ\right)$ for a given partition $\lambda$. As the entries of $\lambda$ grow larger, performing these computations in practice becomes increasingly challenging. However, the general theory provides an algorithmic approach to these calculations.

\subsection*{The computation}\label[section]{theComputation}
Fix a partition $\lambda=(\lambda_1, \lambda_2, \lambda_3, \lambda_4)$. Following \cref{CasiUtili}, we can assume that  $\lambda= (m,s,t,0)$, with $m\geq t+s$. Recall that 
\[
\mathrm{Ext}^k\left(\Sigma_{m,t,s} \mQ,\Sigma_{m,t,s} \mQ\right) \cong H^k(\Sigma_{m,t,s} \mQ \otimes \Sigma_{m,t,s} \mQ^\vee) \cong H^k\left(\Sigma_{m,t,s} \mQ \otimes \Sigma_{m,m-s,m-t} \mQ\otimes \mathcal{O}_X(-m)\right),
\]
hence, using Littlewood--Richardson's formula (\cref{LRD}), we need to compute the cohomology of irreducible factors of the form
\[
E_\mu:=\Sigma_\mu \mQ\otimes \mathcal{O}_X(-m).
\]
To do so, we twist the Koszul complex \eqref{Koszul} by $\Sigma_\mu \Tilde{\mQ}\otimes \mathcal{O}_\mathrm{Gr}(-m)$, obtaining a long exact sequence 
\begin{align*}
0&\to\mathcal{O}_\mathrm{\Gr(6,10)}(-10-m)\to\bigwedge^{19}\bigwedge^3\Tilde{\mU}\otimes \Sigma_\mu \Tilde{\mQ}\otimes \mathcal{O}_\mathrm{Gr}(-m)\to... \\
...&\to\bigwedge^2\bigwedge^3\Tilde{\mU}\otimes \Sigma_\mu \Tilde{\mQ}\otimes \mathcal{O}_\mathrm{Gr}(-m)\to\bigwedge^3\Tilde{\mU} \otimes \Sigma_\mu \Tilde{\mQ}\otimes \mathcal{O}_\mathrm{Gr}(-m)\to \Sigma_\mu \Tilde{\mQ}\otimes \mathcal{O}_\mathrm{Gr}(-m)\to E_\mu\to0.
\end{align*}

As explained in \cref{rmk:dWtoSchur}, each exterior power decomposes into irreducible pieces. For instance, in degree $3$ we will have
\begin{align*}
\Schur_\mu \Tilde{\mQ}\otimes \Schur_{3+m,3+m,1+m,1+m,1+m,m}\Tilde{\mU}&\oplus\Schur_\mu \Tilde{\mQ}\otimes \Schur_{3+m,2+m,2+m,2+m,m,m}\Tilde{\mU}\oplus\\
&\oplus \Schur_\mu \Tilde{\mQ}\otimes \Schur_{2+m,2+m,2+m,1+m,1+m,1+m}\Tilde{\mU}.
\end{align*}
For each irreducible factor appearing in the twisted Koszul complex, it is possible to compute its cohomology using Bott's Theorem.

Denote by $(d^p_0,d^p_1,...,d^p_{24})$ the tuple containing the dimensions of the cohomology vector spaces of 
\[
\bigwedge^{p}\bigwedge^3\Tilde{\mU}\otimes \Sigma_\mu \Tilde{\mQ}\otimes \mathcal{O}_\mathrm{Gr}(-m).
\]

Recall that, for an irreducible vector bundle, at most one cohomology group is non-zero. In particular, lots of the $d^p_i$ could be zero.
To compute the cohomology of $E_\mu$ we split the long exact sequence into short exact sequences:
\[
0\to \bigwedge^{20}\bigwedge^3\Tilde{\mU}\otimes \Sigma_\mu \Tilde{\mQ}\otimes \mathcal{O}_\mathrm{Gr}(-m)\to \bigwedge^{19}\bigwedge^3\Tilde{\mU}\otimes \Sigma_\mu \Tilde{\mQ}\otimes \mathcal{O}_\mathrm{Gr}(-m)\to K_1\to 0
\]
\[
0\to K_1\to \bigwedge^{18}\bigwedge^3\Tilde{\mU}\otimes \Sigma_\mu \Tilde{\mQ}\otimes \mathcal{O}_\mathrm{Gr}(-m)\to K_2\to 0
\]
\[
...
\]
\[
0\to K_{p-1}\to \bigwedge^{20 -p}\bigwedge^3\Tilde{\mU}\otimes \Sigma_\mu \Tilde{\mQ}\otimes \mathcal{O}_\mathrm{Gr}(-m)\to K_p\to 0
\]
\[
...
\]
\[
0\to K_{18}\to \bigwedge^3\Tilde{\mU}\otimes \Sigma_\mu \Tilde{\mQ}\otimes \mathcal{O}_\mathrm{Gr}(-m)\to K_{19}\to 0
\]
\[
0\to K_{19}\to \Sigma_\mu \Tilde{\mQ}\otimes \mathcal{O}_\mathrm{Gr}(-m)\to E_\mu\to 0.
\]

\begin{rmk}\label[remark]{rmk_sol_SES}
   In principle, one can compute the cohomology of $K_1$ using the first short exact sequence and use this information to compute the cohomology of $K_2$ and so on. Eventually, from the last short exact sequence one obtains the cohomology of $E_\mu$, as wanted. In practice, this is possible only for low values of $m$, $t$ and $s$, where the cohomological maps must be trivial.
\end{rmk}
   
\begin{definition}\label[definition]{def_indeterminate}
    We will call a Schur functor $\Sigma_\lambda\mQ$ \emph{indeterminate} if there is a direct summand $E_\mu$ in the decomposition of $\mathcal{E}_\lambda$ with the following property: for at least one value of $p=1,...,20$, there exists a $q> 0$ such that  
    \[
    H^q\left(\mathrm{Gr}(6,10),K_{p-1}\right)\neq 0\neq H^q\left(\mathrm{Gr}(6,10),\bigwedge^{20 -p}\bigwedge^3\Tilde{\mU}\otimes \Sigma_\mu \Tilde{\mQ}\otimes \mathcal{O}_\mathrm{Gr}(-m)\right).
    \]
\end{definition}

Note that, for a bundle $E_\mu$ with the above property one gets a short exact sequence whose induced maps in cohomology of the form 
\[
H^q\left(\mathrm{Gr}(6,10),K_{p-1}\right)\longrightarrow H^q\left(\mathrm{Gr}(6,10),\bigwedge^{20 -p}\bigwedge^3\Tilde{\mU}\otimes \Sigma_\mu \Tilde{\mQ}\otimes \mathcal{O}_\mathrm{Gr}(-m)\right)
\]
are not necessarily trivial. This phenomenon stops the algorithmic process described in \cref{rmk_sol_SES} and we are not able compute the $\mathrm{Ext}$-group of Schur functors $\Schur_{\lambda}\mQ$ if it is indeterminate in the above sense.

\medskip
If $\Schur_{(m,t,s,0)}\mQ$ is not indeterminate, summing up (with multiplicities) all the cohomology of the irreducible factors of the form $E_\mu$ appearing in the Littlewood--Richardson's formula, one obtains $H^*(\Sigma_{m,t,s} \mQ \otimes \Sigma_{m,t,s} \mQ^\vee)$. Following this process, we get the following
\begin{proposition}\label[proposition]{prop:Tab}
    The dimensions of the $\mathrm{Ext}$-group of Schur functors $\Schur_{(m,t,s,0)}\mQ$ with $m < 5$ are the ones stated in \cref{MainTable}. Moreover, any Schur functor of $\mQ$ with $m \geq 5$ is indeterminate.
\end{proposition} 
\begin{proof}
    The first part is obtained with a case-by-case analysis following the previous algorithm. In \cref{WorkedExample}, we treat one example in detail. 
    
    For the second part, we need the next lemma (\cref{GiardinoDellaFelicita}), which tells us that if a factor appears in the Littlewood--Richardson's decomposition of the endomorphisms bundle of a given Schur functor, it keeps appearing as we increase the entries of the partition. We use this lemma to ensure that we did not miss any ``computable'' Schur functors. Take any $\lambda=(m,t,s,0)$ with $m\geq t+s$, we have the following cases:
    \begin{itemize}
        \item \underline{Case $m-t \geq 5$}\\ The indeterminacy follows from the one of $(5,0,0,0)$ by adding $s$-times $(1,1,1,0)$, $(t-s)$-times $(1,1,0,0)$ and $(m-t-5)$-times $(1,0,0,0)$;
        \item \underline{Case $m-t=4$}\\ If $t-s \geq 1$ then the indeterminacy follows from the one of $(5,1,0,0)$, otherwise we check a finite number of cases with $t=s$, $s\leq 4$;
        \item \underline{Case $m-t=3$}\\ If $t-s \geq 2$ then the indeterminacy follows from the one of $(5,2,0,0)$, otherwise we check a finite number of cases with $t=s+1$, $s\leq 3$ and $t=s$, $s\leq 3$;
        \item\underline{Case $m-t=2$}\\ If $t-s \geq 2$ then the indeterminacy follows from the one of $(4,2,0,0)$, otherwise we check a finite number of cases with $t=s+1$, $s\leq 2$ and $t=s$, $s\leq 2$;
        \item\underline{Case $m-t=1$}\\ If $t-s \geq 3$ then the indeterminacy follows from the one of $(4,3,0,0)$, otherwise we check a finite number of cases with $t=s+2$, $s\leq 1$ and $t=s+1$, $s\leq 1$ and $t=s$, $s\leq 1$;
        \item\underline{Case $m-t=0$}\\ In this case $s$ must be zero and the indeterminacy follows from the one of $(3,3,0,0)$ (if $m=t\geq 3$).
    \end{itemize}
\end{proof}

\begin{lemma}\label[lemma]{GiardinoDellaFelicita}
        If $\mu=(\mu_1,\mu_2,\mu_3,\mu_4)$ is a partition appearing in the Littlewood--Richardson decomposition of $\mathcal{E}nd\left(\Schur_{(m,t,s,0)}\mQ\right)$, then $(\mu_1+1,\mu_2+1,\mu_3+1,\mu_4+1)$ appears in the decomposition of $$\mathcal{E}nd\left(\Schur_{(m+1,t,s,0)}\mQ\right), \mathcal{E}nd\left(\Schur_{(m+1,t+1,s,0)}\mQ\right) \text{ and } \mathcal{E}nd\left(\Schur_{(m+1,t+1,s+1,0)}\mQ\right).$$
    \end{lemma}
    \begin{proof}
 This is a generalization of \cite[Proposition 8.3]{FatOno}, where it is proved for $\mathcal{E}nd \left(\Schur_{(m+1,t,s,0)}\mQ\right)$. The other two cases are obtained with the same ideas, following the Littlewood--Richardson's rule.
\end{proof}

\begin{rmk}\label[rmk]{SpectralSeq}
    Notice that studying the Koszul resolution by decomposing it into short exact sequences is equivalent to analyzing the spectral sequence associated with the Koszul complex. 

More precisely, let $A$ be a vector bundle on $\Gr(6,10)$, and suppose we want to compute the cohomology groups $H^p(X, A\vert_X)$. Consider the Koszul complex twisted by $A$:
\[
K_A^{\bullet}:\quad 0 \to A\otimes \bigwedge^{20}\bigwedge^3 \widetilde{\mU} \to \dots \to A\otimes \bigwedge^3 \widetilde{\mU} \to A \to A\vert_X \to 0.
\]
Then
\[
K_A^q = A\otimes \bigwedge^q \bigwedge^3 \widetilde{\mU}, 
\qquad 
E_1^{p,q} = H^p\big(\Gr(6,10), K_A^q\big).
\]
By \cite[Remark 2.67]{Huy06}, this spectral sequence satisfies
\[
E_1^{p,q} = H^p\big(\Gr(6,10), K_A^q\big) \;\Longrightarrow\; 
H^{p+q}\big(\Gr(6,10), A\otimes \of_X\big) = H^{p+q}(X, A\vert_X).
\]

We choose to present the argument using short exact sequences to make the computational differences with respect to \cite{FatOno} more explicit.
\end{rmk}

\subsection{An example of indeterminacy}\label[section]{indeterminacy}
In some cases, when solving the exact sequences in the above process, the cohomology of the irreducible factors forces all the cohomology maps to be trivial. This happens for most of the cases in \cref{MainTable}. However, as the entries of the partition grow, this is no longer the case. Here we present in details an example of this situation for $\lambda=(3,2,1,0)$.

Among the factors of the Littlewood--Richardson decomposition of  $\mathcal{E}nd\left(\Schur_\lambda \mQ\right)$ appears only one which contributes to the indeterminacy, namely
    $$\Schur_{5,5,2,0|3,3,3,3,3,3}.$$ We follow the argument in \cref{theComputation}. By Borel--Weil--Bott’s theorem, the only non-vanishing cohomologies in the Koszul complex are:
    
    \begin{itemize}
        \item $H^{15}\left(\bigwedge^{13}\bigwedge^3\widetilde{\mU}\otimes\Schur_{5,5,2,0}\widetilde{\mQ}\otimes \mathcal{O}_\mathrm{Gr}(-3)\right) = \bigwedge^9V_{10}$;
        \item $H^{14}\left(\bigwedge^{12}\bigwedge^3\widetilde{\mU}\otimes\Schur_{5,5,2,0}\widetilde{\mQ}\otimes \mathcal{O}_\mathrm{Gr}(-3)\right) = \bigwedge^6V_{10}$;
        \item $H^{12}\left(\bigwedge^{11}\bigwedge^3\widetilde{\mU}\otimes\Schur_{5,5,2,0}\widetilde{\mQ}\otimes \mathcal{O}_\mathrm{Gr}(-3)\right) = \mathrm{Sym}^3V_{10}$;
        \item $H^{11}\left(\bigwedge^9\bigwedge^3\widetilde{\mU}\otimes\Schur_{5,5,2,0}\widetilde{\mQ}\otimes \mathcal{O}_\mathrm{Gr}(-3)\right) = \Schur_{2,2,2,2,2,2,2,2,1,0}V_{10}$;
        \item $H^{10}\left(\bigwedge^8\bigwedge^3\widetilde{\mU}\otimes\Schur_{5,5,2,0}\widetilde{\mQ}\otimes \mathcal{O}_\mathrm{Gr}(-3)\right) = \bigwedge^4V_{10}\oplus \Schur_{2,1,1}V_{10}$;
        \item $H^6\left(\bigwedge^4\bigwedge^3\widetilde{\mU}\otimes\Schur_{5,5,2,0}\widetilde{\mQ}\otimes \mathcal{O}_\mathrm{Gr}(-3)\right) = \bigwedge^2V_{10}\oplus \Schur_{2,2,2,1,1,1,1,1,1,0}V_{10}$.
    \end{itemize}
    
The analysis starts with the short exact sequence 
\[
0\to K_6\to \bigwedge^{13}\bigwedge^3\widetilde{\mU}\otimes\Schur_{5,5,2,0}\widetilde{\mQ}\otimes \mathcal{O}_\mathrm{Gr}(-3)\to K_7\to 0,
\]
which tells us that  $H^{15}\left(K_7\right) = \bigwedge^9V_{10}$. Then we have 
\[
0\to K_7\to \bigwedge^{12}\bigwedge^3\widetilde{\mU}\otimes\Schur_{5,5,2,0}\widetilde{\mQ}\otimes \mathcal{O}_\mathrm{Gr}(-3)\to K_8\to 0,
\]
which tells us that  $H^{14}\left(K_8\right) = \bigwedge^9V_{10}\oplus\bigwedge^6V_{10}$. Then we have 
\[
0\to K_8\to \bigwedge^{11}\bigwedge^3\widetilde{\mU}\otimes\Schur_{5,5,2,0}\widetilde{\mQ}\otimes \mathcal{O}_\mathrm{Gr}(-3)\to K_9\to 0,
\]
from which we get that $H^{12}\left(K_{10}\right) = H^{13}\left(K_9\right) = \bigwedge^9V_{10}\oplus\bigwedge^6V_{10}$ and $H^{11}\left(K_{10}\right) = H^{12}\left(K_9\right) = \mathrm{Sym}^3V_{10}$. 

Here comes the first non-trivial short exact sequence
\[
0\to K_{10}\to \bigwedge^9\bigwedge^3\widetilde{\mU}\otimes\Schur_{5,5,2,0}\widetilde{\mQ}\otimes \mathcal{O}_\mathrm{Gr}(-3)\to K_{11}\to 0.
\]
In cohomology, we have
\[
0\to H^{10}(K_{11}) \to \mathrm{Sym}^3V_{10} \xrightarrow{\varphi} \Schur_{2,2,2,2,2,2,2,2,1,0}V_{10} \to H^{11}(K_{11}) \to \bigwedge^9V_{10}\oplus\bigwedge^6V_{10} \to 0.
\]
Using Weyl's formula, we compute $\mathrm{dim}\left(\mathrm{Sym}^3V_{10}\right)=220$ and $\mathrm{dim}\left(\Schur_{2,2,2,2,2,2,2,2,1,0}V_{10}\right)=330$. Therefore, the rank of the map $\varphi$ is not determined by the dimension of the cohomologies. 

\begin{rmk}
    If one can show that $\varphi$ is of maximal rank, it would led to $H^{10}(K_{11})=0$ and $$H^{11}(K_{11})\cong \bigwedge^9V_{10}\oplus\bigwedge^6V_{10} \oplus \faktor{\Schur_{2,2,2,2,2,2,2,2,1,0}V_{10}}{\mathrm{Sym}^3V_{10}}.$$
The remaining short exact sequences are all trivial, thus we would get that the only non-vanishing cohomology of $\Schur_{5,5,2,0|3,3,3,3,3,3} $ is $H^2 \left(\Schur_{5,5,2,0|3,3,3,3,3,3}\right)$, given by
\[
\bigwedge^9V_{10}\oplus\bigwedge^6V_{10} \oplus \faktor{\Schur_{2,2,2,2,2,2,2,2,1,0}V_{10}}{\mathrm{Sym}^3V_{10}} \oplus \bigwedge^4V_{10}\oplus \Schur_{2,1,1}V_{10} \oplus \bigwedge^2V_{10}\oplus \Schur_{2,2,2,1,1,1,1,1,1,0}V_{10},
\] 
whose dimension is $2730$. If that is the case, and after computing the $\mathrm{Ext}$ cohomology of all the other factors of the Littlewood--Richardson decomposition of $\mathcal{E}nd\left(\Schur_\lambda \mQ\right)$, one obtains that 
\[
\mathrm{ext}^0(\Sigma_\lambda \mQ,\Sigma_\lambda \mQ) = 1, \quad \mathrm{ext}^1(\Sigma_\lambda \mQ,\Sigma_\lambda \mQ) = 40, \quad \mathrm{ext}^2(\Sigma_\lambda \mQ,\Sigma_\lambda \mQ) = 35406.
\]
These dimensions are the same as the ones computed for the correspondent partition on the Beauville--Donagi fourfold in \cite{FatOno}. In fact, in that case the Schur functor $\Schur_{(3,2,1,0)}(\mQ)$ is not indeterminate.
\end{rmk}

\subsection{Combinatorial results}
Here we deal with some combinatorial aspects of the cohomological dimensions of the sheaf of endomorphisms of Schur functors. The goal is to follow the Littlewood--Richardson's rule to predict which Schur functors appear in the decomposition of $\mathcal{E}nd\left(\Schur_\lambda\mQ\right)$. The main result is \cref{multiplicitiesLDR}.
\begin{lemma}\label[lemma]{LRDComputation}
Let $\lambda = (m,t,s,0)$ be a partition with $m\geq t+s$, then, among the irreducible factors of $$\Schur_{(m,t,s,0)}\otimes \Schur_{(m,m-s,m-t,0)},$$ we have:
    \begin{itemize}[(1)]
    \item[(1)] $\Schur_{(m,m,m,m)}$ with multiplicity $1$;
    \item[(2)] $\Schur_{(m+1,m+1,m-1,m-1)}$ with multiplicity $
    \begin{cases}
        0 & \text{ if ~$t=s=0$}\\
        2 & \text{ if ~$m>t>s>0$}\\
        1 & \text{ otherwise}
    \end{cases}
    $
    \item[(3)] $\Schur_{(m+2,m,m,m-2)}$ with multiplicity $
    \begin{cases}
        0 & \text{ if ~$m=1$}\\
        \geq 1 & \text{ otherwise}
    \end{cases}
    $
\end{itemize}
\end{lemma}
\begin{proof}
 This is an application of Littlewood--Richardson's formula, recall the three rules (a), (b) and (c) of \cref{LRD}. Note that, in order to follow (c), it is more convenient to fill the diagrams from right to left, top to bottom. In the following, we call \textit{admissible fillings} any filling of a Young diagram which respects the Littlewood--Richardson's rules.
    \begin{multicols}{2}
        For (1), we need to show that there is exactly one way to fill the diagram on the right with content $(m,m-s,m-t,0)$.
        
        \[
    \ytableausetup{centertableaux}
    \ytableaushort
    {\none\none\none \none\none \dots,\none\none \dots, \dots }
    * {7,7,7,7}
    * [*(gray)]{7,4,1,0}
    \]
    \end{multicols}
    By (c) we have to start with a $1$ and the columns have to be strictly increasing for (b); since by (a) every row must be weakly increasing, we get
    \[
    \ytableausetup{nosmalltableaux}
    \ytableaushort
    {\none\none\none\none\none\dots,\none\none\dots\none 1\dots 1, \dots }
    * {7,7,7,7}
    * [*(gray)]{7,4,1,0}
    ~~~~~~~~~~~~~~~~~~~~~~~~~~~~~~~~
    \ytableaushort
    {\none\none\none\none\none\dots,\none\none\dots\none 1\dots 1, 
    \dots\none\none\none 2\dots 2}
    * {7,7,7,7}
    * [*(gray)]{7,4,1,0}
    \]
    Now, the number of $2$ equals the number of $1$, and hence we complete the third and fourth rows as
    \[
    \ytableausetup{centertableaux}
    \ytableaushort
    {\none\none\none\none\none\dots,\none\none\dots\none 1\dots 1, 
    \dots 1 \dots 1 2\dots 2}
    * {7,7,7,7}
    * [*(gray)]{7,4,1,0}
    ~~~~~~~~~~~~~~~~~~~~~~~~~~~~~~~~
    \ytableausetup{centertableaux}
    \ytableaushort
    {\none\none\none\none\none\dots,\none\none\dots\none 1\dots 1, 
    \dots 1 \dots 1 2\dots 2,
    1 2 \dots 2 3 \dots 3}
    * {7,7,7,7}
    * [*(gray)]{7,4,1,0}
    \]
    One can check that the content of this semi--standard Young tableaux is $(m,m-s,m-t,0)$. Since we did not made any choice during the filling process, the multiplicity of $\Schur_{(m,m,m,m)}$ in $\Schur_{(m,t,s,0)}\otimes \Schur_{(m,m-s,m-t,0)}$ is $1$.
    
    In order to prove (2), we consider five separate cases.
    \begin{itemize}
        \item \underline{Case 1.} 
        \begin{multicols}{2}
            Assume that $t=s=0$.\\
            We want to show that there is no admissible filling for the diagram on the right.\\
            \[
    \ytableausetup{centertableaux}
    \ytableaushort
    {\none\dots,\none\dots,\none\dots,\none\dots}
    * {5,5,3,3}
    * [*(gray)]{4,0,0,0}
    \]
        \end{multicols} 
    
    \begin{multicols}{2}
        By (c) we have a $1$ in the first row and then a $2$ under it by (b). Hence we get the tableaux on the right, which is not admissible since we have $m+1$ times the number $1$.\\
    \[
    \ytableausetup{centertableaux}
    \ytableaushort
    {\none\dots\none\none 1, 1 \dots 1 1 2,\none\dots,\none\dots}
    * {5,5,3,3}
    * [*(gray)]{4,0,0,0}
    \]
    \end{multicols}
        \item \underline{Case 2.} 
        \begin{multicols}{2}
            Assume that $m=t$, then $s=0$\\
           The only admissible filling for the diagram is the one on the right.\\
            \[
    \ytableausetup{centertableaux}
    \ytableaushort
    {\none\dots\none\none1,\none\dots\none\none 2,1\dots1,2\dots2}
    * {5,5,3,3}
    * [*(gray)]{4,4,0,0}
    \]
        \end{multicols}

    \item \underline{Case 3.} 
    \begin{multicols}{2}
        Assume that $t>s=0$.\\
    We want to show that there is exactly one way to fill the diagram on the right.\\
    \[
    \ytableausetup{centertableaux}
    \ytableaushort
    {\none\dots\none\none\none\dots,
    \none\dots\none\none\none\dots,
    \none\dots\none\none\none\dots,
    \none\dots\none\none\none\dots}
    * {10,10,8,8}
    * [*(gray)]{9,4,0,0}
    \]
    \end{multicols}
        We start as before. For the third row, we put the $2s$ under all the $1s$, along with exactly one other $2$: for (c), this is the last admissible $2$ at this step and if we don't put it then we will have too much $1s$ as before.
        \[
    \ytableausetup{centertableaux}
    \ytableaushort
    {\none\dots\none\none\none\dots\none\none\none 1,
    \none\dots\none\none 1 \dots 1 1 1 2,
    \none\dots\none\none\none\dots,
    \none\dots\none\none\none\dots}
    * {10,10,8,8}
    * [*(gray)]{9,4,0,0}
~~~~~~~~~~~~~~~~
    \ytableausetup{centertableaux}
    \ytableaushort
    {\none\dots\none\none\none\dots\none\none\none 1,
    \none\dots\none\none 1 \dots 1 1 1 2,
    1 \dots 1 2 2 \dots 2 2,
    \none\dots\none\none\none\dots}
    * {10,10,8,8}
    * [*(gray)]{9,4,0,0}
    ~~~~~~~~~~~~~~~~
    \ytableausetup{centertableaux}
    \ytableaushort
    {\none\dots\none\none\none\dots\none\none\none 1,
    \none\dots\none\none 1 \dots 1 1 1 2,
    1 \dots 1 2 2 \dots 2 2,
    2 \dots 2 3 3 \dots 3 3}
    * {10,10,8,8}
    * [*(gray)]{9,4,0,0}
    \]
    \item \underline{Case 4.}
    \begin{multicols}{2}
        Assume that $t=s>0$. We want to show that there is exactly one way to fill the diagram on the right.\\
    \[
    \ytableausetup{centertableaux}
    \ytableaushort
    {\none\dots\none\none\none\none\dots,
    \none\dots\none\none\none\none\dots,
    \none\dots\none\none\none\none\dots,
    \none\dots\none\none\none\none\dots}
    * {10,10,8,8}
    * [*(gray)]{9,4,4,0}
    \]
    \end{multicols}
    
    As before, we get the following filling.
    \[
    \ytableausetup{centertableaux}
    \ytableaushort
    {\none\dots\none\none\none\none\dots\none\none 1,
    \none\dots\none\none 1 1 \dots 1 1 2,
    \none\dots\none\none\none\none\dots,
    \none\dots\none\none\none\none\dots}
    * {10,10,8,8}
    * [*(gray)]{9,4,4,0}
    ~~~~~~~~~~~~~~~~
    \ytableausetup{centertableaux}
    \ytableaushort
    {\none\dots\none\none\none\none\dots\none\none 1,
    \none\dots\none\none 1 1 \dots 1 1 2,
    \none\dots\none\none 2 2 \dots 2,
    \none\dots\none\none\none\none\dots}
    * {10,10,8,8}
    * [*(gray)]{9,4,4,0}
    ~~~~~~~~~~~~~~~~
    \ytableausetup{centertableaux}
    \ytableaushort
    {\none\dots\none\none\none\none\dots\none\none 1,
    \none\dots\none\none 1 1 \dots 1 1 2,
    \none\dots\none\none 2 2 \dots 2,
    1 \dots 1 3 3 3 \dots3}
    * {10,10,8,8}
    * [*(gray)]{9,4,4,0}
    \]
    \item \underline{Case 5.}
    \begin{multicols}{2}
        Assume that $t>s>0$. There are two possible ways to fill the diagram on the right.\\
    
    \[
    \ytableausetup{centertableaux}
    \ytableaushort
    {\none\dots\none\none\none\dots\none\none\none\dots,
    \none\dots\none\none\none\dots\none\none\none\dots,
    \none\dots\none\none\none\dots\none\none\none\dots,
    \none\dots\none\none\none\dots\none\none\none\dots}
    * {13,13,11,11}
    * [*(gray)]{12,8,4,0}
    \]
    \end{multicols}

   \[
    \ytableausetup{centertableaux}
    \ytableaushort
    {\none\dots\none\none\none\dots\none\none\none\dots\none\none 1,
    \none\dots\none\none\none\dots\none\none 1 \dots 1 1 2,
    \none\dots\none\none 1 \dots 1 1 2 \dots 2,
    1 \dots 1 2 2 \dots 2 3 3 \dots 3}
    * {13,13,11,11}
    * [*(gray)]{12,8,4,0}
    ~~~~~~~~~~~~~~~~~~~~~~~~~~~~~~~~
    \ytableausetup{centertableaux}
    \ytableaushort
    {\none\dots\none\none\none\dots\none\none\none\dots\none\none 1,
    \none\dots\none\none\none\dots\none\none 1 \dots 1 1 2,
    \none\dots\none\none 1 \dots 1 2 2 \dots 2,
    1 \dots 1 1 2 \dots 2 3 3 \dots 3}
    * {13,13,11,11}
    * [*(gray)]{12,8,4,0}
    \]
    \end{itemize}
\vspace{1cm}
    It remains to prove $(3)$. Note that if $m=1$ then we have only the cases $\lambda=(1,0,0,0)$, $(1,1,0,0)$ and the result trivially holds. Assume that $m>1$. We divided the proof into three cases.
    \begin{itemize}
        \item \underline{Case 1.} Assume that $m-t=0$ then we can fill the diagram (at least) in the following way
        \[
    \ytableausetup{centertableaux}
    \ytableaushort
    {\none\dots\none\none\none 1 1,
    \none\dots\none\none\none,
    1 \dots 1 2 2,
    2 \dots 2}
    * {7,5,5,3}
    * [*(gray)]{5,5,0,0}
    \]
    \item \underline{Case 2.} Assume that $m-t=1$ then we can fill the diagram (at least) in the following ways (respectively if $s =1 $ and if $s=0$).
    \[
    \ytableausetup{centertableaux}
    \ytableaushort
    {\none\none\dots\none\none\none 1 1,
    \none\none\dots\none\none 2,
    \none 1 \dots 1 2 3,
    1 2 \dots 2}
    * {8,6,6,4}
    * [*(gray)]{6,5,1,0}
    ~~~~~~~~~~~~~~~~~~~~~~~~~~~~~~~~
    \ytableausetup{centertableaux}
    \ytableaushort
    {\none\dots\none\none\none 1 1,
    \none\dots\none\none 2,
    1 \dots 1 2 3,
    2 \dots 2}
    * {7,5,5,3}
    * [*(gray)]{5,4,0,0}
    \]
    \item \underline{Case 3.} Assume that $m-t\geq 2$ then we can fill the diagram (at least) in the following way
        \[
    \ytableausetup{centertableaux}
    \ytableaushort
    {\none\dots\none\none\dots\none\none\dots\none\none\none 1 1,
    \none\dots\none\none\dots\none 1 \dots 1 2 2 ,
    \none\dots\none 1 \dots 1 2 \dots 2 3 3,
    1 \dots 1 2 \dots 2 3 \dots 3}
    * {13,11,11,9}
    * [*(gray)]{11,6,3,0}
    \]
    \end{itemize}
    \vspace{0.5cm}
    Note that there are other admissible fillings besides these so that the multiplicity can be bigger than $1$.
\end{proof}

\begin{corollary}\label[corollary]{multiplicitiesLDR}
    Let $\lambda = (m,t,s,0)$ with $m\geq t+s$. Among the factors of the Littlewood--Richardson decomposition of $\Sigma_{\lambda}\widetilde{\mQ}\otimes \Sigma_{\lambda}\widetilde{\mQ}^\vee$ we have
    \begin{itemize}
    \item $\mathcal{O}_\mathrm{Gr}$ with multiplicity $1$;
    \item $\Schur_{(2,2,0,0)}\widetilde{\mQ}\otimes\mathcal{O}_\mathrm{Gr}(-1)$ with multiplicity $
    \begin{cases}
        0 & \text{ if ~$t=s=0$}\\
        2 & \text{ if ~$m>t>s>0$}\\
        1 & \text{ otherwise}
    \end{cases}
    $
    \item $\Schur_{(4,2,2,0)}\widetilde{\mQ}\otimes\mathcal{O}_\mathrm{Gr}(-2)$ with multiplicity $
    \begin{cases}
        0 & \text{ if ~$m=1$}\\
        \geq 1 & \text{ otherwise}
    \end{cases}
    $
\end{itemize}

\end{corollary}
\begin{proof}
    The result follows directly from the \cref{LRDComputation}. In fact, we need to twist each irreducible factor by $\mathcal{O}_\mathrm{Gr}(-m)$ in order to get an irreducible piece of $\Sigma_\lambda \widetilde{\mQ}\otimes \Sigma_\lambda \widetilde{\mQ}^\vee$. 
\end{proof}

\begin{rmk}
    Following the steps described in \cref{theComputation}, we have that 
    \[
    H^1\left(\Schur_{(2,2,0,0)}\mQ\otimes\mathcal{O}_X(-1)\right) \cong V_{10}\oplus V_{10}^\vee \ \text{ and }\ H^2\left(\Schur_{(2,2,0,0)}\mQ\otimes\mathcal{O}_X(-1)\right) \cong \mathbb{C}.
    \]
    One can also notice that, differently from the Beauville--Donagi cases, in these cases the $\Ext^1$'s are not irreducible representations. 
\end{rmk}

\subsection{Stability and the notion of (projectively) hyper-holomorphic vector bundles}\label[section]{sec:hyperholo}
In \cite{Ver96}, Verbitsky introduced the notion of hyper-holomorphic connection for vector bundles and of projectively hyper-holomorphic vector bundle on a \hk manifold $X$. The idea behind the definition is to select those holomorphic bundles that are holomorphic for every complex structure induced by the hyperk\"ahler structure on $X$. Moreover, he managed to prove that (the smooth locus of an irreducible component of) the moduli space parametrizing projectively hyper-holomorphic vector bundles with fixed numerical invariants is hyperk\"ahler, see Theorem 3.11 \textit{op. cit.}.

More recently, the study of infinitesimal deformations of hyper-holomorphic vector bundles has been carried on by Meazzini and Onorati in \cite{MeaOno}.

\begin{rmk}\label[rmk]{ModularHyperHolo}
    For $X$ of type $\textrm{K3}^{[2]}$, if $E$ is a stable modular vector bundle on $X$ then $E$ is projectively hyper-holomorphic. This follows from \cite[Proposition 1.2 and Theorem 2.5]{Ver96}, see also \cite[Proposition 5.6]{Bec25}.
\end{rmk}

\cref{multiplicitiesLDR} is a promising indication that the vector bundles $\Sigma_\lambda \mQ$ are stable. In fact, we have the following result.

\begin{proposition}\label[proposition]{rmk:stab}
    Let $X\subseteq \mathrm{Gr}(6,10)$ be a very general Debarre--Voisin \hk manifold, then $\Sigma_\lambda \mQ$ is polystable for every partition $\lambda$.
\end{proposition}
\begin{proof}
    Recall that a vector bundle is polystable if it is isomorphic to a direct sum of stable bundles. By \cite[Proposition 8.4]{OG22}, $\mQ$ is stable and by \cite[Theorem 3.2.11]{HL} we have that $\mQ^{\otimes k}$ is polystable for every $k\geq 1$. Since every Schur functor $\Sigma_\lambda \mQ$ is a direct summand of $\mQ^{\otimes |\lambda|}$, it follows that it is polystable.
\end{proof}

\begin{rmk}\label[rmk]{ConjStab}
    By the previous result, we have that $$\Sigma_\lambda \mQ\cong \bigoplus_{i =1 }^k \mathcal{F}_i$$ is a sum of stable bundles. Therefore, $\Sigma_\lambda \mQ \otimes \Sigma_\lambda \mQ^\vee$ contains at least $k$ copies of $\mathcal{O}_X$. The combinatorial result \cref{multiplicitiesLDR} tells us that, at the level of $\Gr(2,6)$, one has exactly one copy of $\mathcal{O}_{\mathrm{Gr}(6,10}$.
    
    If the same holds for the restriction to $X$, one gets that $k=1$ and that $\Sigma_\lambda \mQ$ is stable. We remark that our computations suggest that this is the case.

    In particular, if $X\subseteq \mathrm{Gr}(6,10)$ is a very general Debarre--Voisin \hk manifold and if $\Sigma_\lambda \mQ$ is stable, then $\Sigma_\lambda \mQ$ is a stable projectively hyper-holomorphic vector bundle.
\end{rmk}

\section{A worked example}\label[section]{WorkedExample}
This section shows the computations for the modular sheaf $\bigwedge^2\mQ$ with $\mathrm{ext}^1=20$.
Following \cref{notz:endBundle}, consider $\mathcal{E}_{1,1}$, which, by Littlewood--Richardson's decomposition, is equal to $$\Schur_{1,1,0,0}\mQ\otimes\Schur_{1,1,0,0}\mQ\otimes\of_X(-1)= \Sigma_{2,2,0,0}\mQ(-1)\oplus\Sigma_{2,1,1,0}\mQ(-1)\oplus\of_X.$$ In order to compute $h^p(X,\mathcal{E}_{1,1})$, we need to twist the Koszul complex with $\widetilde{\mathcal{E}}_{1,1}$. In this way, we can split it into the three Koszul complexes:
\[
0\to \bigwedge^{20}\bigwedge^3\widetilde{\mU}\otimes\Sigma_{2,2,0,0}\widetilde{\mQ}(-1) \to...\to \bigwedge^3\widetilde{\mU}\otimes\Sigma_{2,2,0,0}\widetilde{\mQ}(-1) \to \Sigma_{2,2,0,0}\widetilde{\mQ}(-1)\to\Sigma_{2,2,0,0}\mQ(-1)\to 0
\]
\[
0\to \bigwedge^{20}\bigwedge^3\widetilde{\mU}\otimes\Sigma_{2,1,1,0}\widetilde{\mQ}(-1) \to...\to \bigwedge^3\widetilde{\mU}\otimes\Sigma_{2,1,1,0}\widetilde{\mQ}(-1) \to \Sigma_{2,1,1,0}\widetilde{\mQ}(-1)\to\Sigma_{2,1,1,0}\mQ(-1)\to 0
\]

\[
0\to \bigwedge^{20}\bigwedge^3\widetilde{\mU}\to...\to\bigwedge^2\bigwedge^3\widetilde{\mU}\otimes\to\bigwedge^3\widetilde{\mU}\to\of_{\Gr(6,10)}\to\of_{X}\to 0,
\]
compute the cohomologies and sum them at the end. By Borel--Weil--Bott’s theorem, the complex twisted by $\Schur_{2,2,0,0}\widetilde{\mQ}(-1)$ has five factors that are not acyclic, namely: 

\begin{align*}
    \Schur_{2,2,0,0|4,4,2,2,2,1} \text{ in } \bigwedge^3\bigwedge^3\widetilde{\mU}\otimes\Schur_{2,2,0,0}\widetilde{\mQ}(-1),\\
    \Schur_{2,2,0,0|8,5,5,3,3,3} \text{ in } \bigwedge^7\bigwedge^3\widetilde{\mU}\otimes\Schur_{2,2,0,0}\widetilde{\mQ}(-1),\\
    \Schur_{2,2,0,0|8,5,5,3,3,3} \text{ in } \bigwedge^7\bigwedge^3\widetilde{\mU}\otimes\Schur_{2,2,0,0}\widetilde{\mQ}(-1),\\
    \Schur_{2,2,0,0|8,8,6,6,4,4}\text{ in } \bigwedge^{10}\bigwedge^3\widetilde{\mU}\otimes\Schur_{2,2,0,0}\widetilde{\mQ}(-1),\\
    \Schur_{2,2,0,0|9,9,9,7,7,4}\text{ in } 
    \bigwedge^{13}\bigwedge^3\widetilde{\mU}\otimes\Schur_{2,2,0,0}\widetilde{\mQ}(-1),\\
    \Schur_{2,2,0,0|11,10,10,10,8,8}\text{ in } \bigwedge^{17}\bigwedge^3\widetilde{\mU}\otimes\Schur_{2,2,0,0}\widetilde{\mQ}(-1).
\end{align*}

Again, by Borel--Weil--Bott’s theorem, all the factors of the complex twisted by $\Sigma_{2,1,1,0}\widetilde{\mQ}(-1)$ are acyclic, and hence do not contribute to the cohomology of $\mathcal{E}_{1,1}$. The third complex has three factors that are not acyclic:

\begin{align*}
    \bigwedge^{20}\bigwedge^3\widetilde{\mU}&=\of_{\Gr(6,10)}(-10),\\
    \Schur_{0,0,0,0|7,7,7,3,3,3}&\text{ in } \bigwedge^{10}\bigwedge^3\widetilde{\mU}\text{ and } \of_{\Gr(6,10)}.
\end{align*}

Putting all together, we obtain: 
\begin{itemize}
\item $H^0(X,\mathcal{E}_{1,1})=H^0\left(\Gr(6,10),\of_{\Gr(6,10)}\right)=\mathbb{C}$;
\item $H^1(X,\mathcal{E}_{1,1})=H^4\left(\Gr(6,10),\Schur_{2,2,0,0|4,4,2,2,2,1}\right)\oplus H^8\left(\Gr(6,10),\Schur_{2,2,0,0|8,5,5,3,3,3}\right)=V_{10}^\vee\oplus V_{10}$;
\item $H^2(X,\mathcal{E}_{1,1})=H^{12}\left(\Gr(6,10),\Schur_{2,2,0,0|8,8,6,6,4,4}\right)\oplus H^{12}\left(\Gr(6,10),\Schur_{0,0,0,0|7,7,7,3,3,3}\right)=\mathbb{C}\oplus\mathbb{C}$;
\item $H^3(X,\mathcal{E}_{1,1})=H^{16}\left(\Gr(6,10),\Schur_{2,2,0,0|9,9,9,7,7,4}\right)\oplus H^{20}\left(\Gr(6,10),\Schur_{2,2,0,0|11,10,10,10,8,8}\right)=V_{10}\oplus V_{10}^\vee$;
\item $H^4(X,\mathcal{E}_{1,1})=H^{24}\left(\Gr(6,10),\of_{\Gr(6,10)}(-10)\right)=\mathbb{C}$.
\end{itemize}
From \cref{SchurIsModular}, we get that $\bigwedge^2\mQ$ is modular, hence we obtained an example of a modular sheaf on the Debarre--Voisin fourfold with $\mathrm{Ext}^1$ of dimension $20$, which was announced in \cite{Fat24}, but without a proof. Note that here $\Ext^2=\mathbb{C}^2$, with $\Ext_0^2=\mathbb{C}$.

In \cite[Section 3.3]{OG26}, O'Grady described the vector bundle $\bigwedge^2\mQ$ on the Beauville--Donagi fourfold in terms of deformations of a modular sheaf on a Hilbert scheme of pairs of points on a $\textrm{K3}$ surface. The same techniques can be adapted also for the Debarre--Voisin case. 

Let $(S, D)$ be a polarized K3 surface with $D\cdot D = 6$ and let $\mathcal{F}$ be a stable spherical vector bundle on $S$ with Mukai vector
\[
v(\mathcal{F}):=(rk(\mathcal{F}),\ch_1(\mathcal{F}),rk(\mathcal{F})+\ch_2(\mathcal{F})) = (2,D,2).
\]
The following is the analogue of \cite[Proposition 3.9]{OG26} for the Debarre--Voisin fourfold.

\begin{proposition}
    Let $S$ be a K3 surface as before. Assume in addition the technical conditions of \cite[Lemma 3.8]{OG26}. If $(X,h)$ is a general deformation of $(S^{[2]},2D-\delta)$, i.e. a general Debarre--Voisin fourfold, then the pair $\left(S^{[2]},\mathcal{G}\left(\bigwedge^2\mathcal{F},\mathrm{Sym}^2\mathcal{F}\right)\right)$ deforms to the pair $\left(X,\bigwedge^2\mQ\right)$. Here $\mathcal{G}$ is the sheaf defined in \cite[Definition 2.1]{OG26}.
\end{proposition}
\begin{proof}
    The proof \cite[Proposition 3.8]{OG26} can be adapt to this case after noticing that the bundle $\mathcal{F}[2]^+$ deforms to $\mQ$ over $X$. This follows from the stability of $\mQ$ and the main result in \cite{OG22}. In fact, this is contained in the proof of \cite[Proposition 8.7]{OG22}. The rest of the proof works verbatim as it only involves properties of the sheaf $\mathcal{F}[2]^+$ on $S^{[2]}$.
\end{proof}

\begin{rmk}
    This correspondence could be a starting point to show that the Schur functors of $\mQ$ on the Beauville--Donagi fourfold deforms to Schur functors of $\mQ$ on the Debarre--Voisin fourfold, possibly passing through some modular sheaves on the Hilbert scheme of points on a $\textrm{K3}$ surface.
\end{rmk}

\section{About the symmetric power}\label[section]{AboutSym}
In this section, we will restrict ourselves to the case of $\mathrm{Sym}^m(Q)$, which is a particular Schur functor given by $\lambda = (m,0,0,0)$. We perform some of the computations of the previous sections in this specific case to make them clearer and concrete. 

We start by computing the Chern classes of $\mathrm{Sym}^m(Q)$.
\begin{lemma}\label[lemma]{ChernSym}
Let $X\subseteq \mathrm{Gr}(6,10)$ a very general Debarre--Voisin fourfold, then:
\begin{enumerate}
    \item $\mathrm{rk}(\mathrm{Sym}^mQ)=\binom{m+3}{3}\frac{(m+3)(m+2)(m+1)}{6}=:r_m$;
    \item $\ch(\mathrm{Sym}^mQ)= r_m + \left(\frac{m}{4}r_m\ch_1(Q)\right)+ \left( \frac{m^2-m}{40}r_m \ch_1(Q)^2 + \frac{m^2+4m}{20}r_m\ch_2(Q)\right) + \\ + \left( -\frac{2m^3-3m^2}{4}r_m \ch_3(Q)\right) + \left( -\frac{10m^4-30m^3+21m^2-6m}{20} r_m \ch_4(Q)\right)$;
    \item $\Delta(\mathrm{Sym}^mQ) = \frac{m^2+4m}{80}r_m^2 \Delta(Q)=\frac{m^2+4m}{80}r_m^2 \mathsf{c}_2(X)$;
    \item $\Xi(\mathrm{Sym}^mQ)=\frac{m(m+4)(9m^2+36m-5)}{110400}r_m^2 \mathsf{c}_2(X)^2$;
    \item $\int_X\Xi(\mathrm{Sym}^mQ)=\frac{3m(m+4)(9m^2+36m-5)}{400}r_m^2$.
    \item $\chi(\Sym^m\mQ,\Sym^m\mQ)=3\left(\frac{3m^2+12m-20}{20}\right)^2r^2_m$.
    
\end{enumerate}
\end{lemma}
\begin{proof} We proceed in order:
\begin{enumerate}
    \item it follows from the dimension of the $m$-symmetric power of a $4$ dimensional vector space;
    
    \item it follows from a more general result in \cite{Svr93}. One needs to apply to our case \cite[Theorem 4.8]{Svr93};
    
    \item it follows directly from the definitions;

    \item it follows from the definitions too;
    
    \item it is a combination of $(1)$ and $(4)$ of \cref{toddIHS};
    
    \item it is an application of Hirzebruch--Riemann--Roch theorem:
    \[
    \chi(\Sym^m\mQ,\Sym^m\mQ)= \chi(\Sym^m\mQ\otimes \Sym^m\mQ^\vee)= \int_X\ch\left(\mathcal{E}nd(\Sym^m\mQ)\right) \mathrm{td}_X
    \]
    Recall that $\ch\left(\mathcal{E}nd(\Sym^m\mQ)\right)= r_m^2 + 0 - \Delta(\Sym^m\mQ) + 0 + \Xi(\mathrm{Sym}^mQ)$. Thus,
    \begin{align*}
     &\chi(\Sym^m\mQ,\Sym^m\mQ)=\int_X\left(3r_m^2\mathsf{p}-\frac{1}{12}\Delta(\Sym^m\mQ)\mathsf{c}_2(X) + \Xi(\mathrm{Sym}^mQ)\right) = \\
     &= 3r_m^2-\frac{828}{12}\frac{m^2+4m}{80}r_m^2+ \frac{3m(m+4)(9m^2+36m-5)}{400}r_m^2 = 3\left(\frac{3m^2+12m-20}{20}\right)^2r^2_m.
    \end{align*}
\end{enumerate}
\end{proof}

Let us compute the cohomology of the $\Ext$-group of $\Sym^mQ$. Recall that we denoted the sheaf of endomorphisms $\mathcal{E}nd\left(\mathrm{Sym}^m\mQ,\mathrm{Sym}^m\mQ\right)$ by $\mathcal{E}_m$.
\begin{rmk}\label[rmk]{rmk:sym-skip}
By Pieri's rule, we get the following decomposition:
\[
    \mE_m=\bigoplus_{i=0}^m\Schur_{2m-i,m,m,i}Q\otimes \mathcal{O}_X(-m),
\]

which leads to $\mE_m=\Schur_{2m,m,m,0}\mQ \otimes \mathcal{O}_X(-m) \oplus\mE_{m-1}$. This means that we can care just about $E_m:=\Schur_{2m,m,m,0}\mQ \otimes \mathcal{O}_X(-m) $ for each $m$. 

\end{rmk}

We now state a useful lemma which spares us some computations later.

\begin{lemma}\label[lemma]{lm:dropcomp}
Let be $E_m$ as above. If we twist the Koszul complex \eqref{Koszul} by $\widetilde{E}_m:=\Schur_{2m,m,m,0}\widetilde{\mQ} \otimes \mathcal{O}_{\Gr(6,10)}(-m)$, we have the resolution:
    \[
        0\to \widetilde{E}_m\otimes\bigwedge^{20}\bigwedge^{3}\widetilde{\mU}\to...\to\widetilde{E}_m\otimes\bigwedge^{2}\bigwedge^{3}\widetilde{\mU}\to\widetilde{E}_m\otimes\bigwedge^{3}\widetilde{\mU}\to\widetilde{E}_m\to E_m\to 0. 
    \]\label[eq]{symKoszul}
    
    Denote $\widetilde{E}_m\otimes\bigwedge^{q}\bigwedge^{3}\widetilde{\mU}$ by $E_m^q$, for $0\leq q\leq 10$. Then, among the irreducible factors of $E_m^q=\bigoplus_i\Sigma_{\left(2m,m,m,0|\mu_{m,i}^q\right)}$, the following partitions $\mu_m^q$ are the ones for which $\Schur_{(2m,m,m,0|\mu_{m,i}^q)}$ is not acyclic for all $m$:
    
    \begin{table} 
    \centering
    \begin{adjustbox}{max width=\textwidth}
    \begin{tabular}{|c|c|c|c|c|}
        \hline
          & $\mu^q_{m,1}$ & $\mu^q_{m,2}$ & $\mu^q_{m,3}$ & $\mu^q_{m,4}$ \\ 
        \hline
         $q=0$ & $(m,m,m,m,m,m)$ & - & - & - \\ 
        \hline
         $q=1$ & $(m+1,m+1,m+1,m,m,m)$ & - & - & - \\ 
        \hline
         $q=2$ & $(m+1,m+1,m+1,m+1,m+1,m+1)$ & - & - & - \\ 
        \hline
         $q=6$ & $(m+5,m+5,m+3,m+3,m+1,m+1)$ & $(m+5,m+5,m+2,m+2,m+2,m+2)$ & - & - \\ 
        \hline
         $q=7$ & $(m+6,m+5,m+3,m+3,m+2,m+2)$ & $(m+5,m+5,m+3,m+3,m+3,m+2)$ & - & - \\ 
        \hline
         $q=8$ & $(m+6,m+6,m+3,m+3,m+3,m+3)$ & - & - & - \\ 
        \hline
         $q=9$ & $(m+7,m+6,m+6,m+3,m+3,m+2)$ & $(m+6,m+6,m+6,m+4,m+4,m+1)$ & - & - \\ 
        \hline
         $q=10$ & $(m+8,m+6,m+6,m+4,m+4,m+2)$ & $(m+8,m+6,m+6,m+4,m+3,m+3)$ & $(m+7,m+7,m+7,m+3,m+3,m+3)$ & $(m+7,m+7,m+6,m+4,m+4,m+2)$ \\ 
        \hline
    \end{tabular}
    \end{adjustbox}
    \newline 
    \caption{}
    \label[table]{tab:sym-transposed}
    \end{table}
    Moreover \[E^{10+k}_m=\widetilde E_m\otimes\bigwedge^{10+k}\bigwedge^3\widetilde\mU\cong \widetilde E_m\otimes\bigwedge^{10-k}\bigwedge^3\widetilde\mU\otimes\of(-k)=E_m^{10-k}(-k), \ \textit{for} k\in\Z_{\geq 1},\]thus $E_{m}^{10-k}$ and $E_m^{10+k}$ have the same number of irreducible factors that are not acyclic for all $m$.
\end{lemma}
\begin{proof}
    Let $\mu_m^q$ be a partition appearing in \cref{tab:your_label} twisted by $\of(-m)$. Consider the partition $(2m,m,m,0|\mu_m^q)$, if we sum the partition $\delta=(9,8,7,6,5,4,3,2,1,0)$, then we have $(2m,m,m,0|\mu_m^q)+\delta=(2m+9,m+8,m+7,6|\mu_m^q+\delta')$ with $\delta'=(5,4,3,2,1,0)$. By confronting $(2m+9,m+8,m+7,6)$ with $\mu_m^q+\delta'$ we recollect only the partitions with no repeated entries. By \cref{rmk:dWtoSchur} we can restrict ourselves for $0\leq q\leq 10$.
\end{proof}

Now that we know which factors are involved in the cohomologies we can make $m$ running in $\mathbb{N}$ and recover a result analogous to \cite[Proposition 2.3]{FatOno}:

\begin{proposition}\label[proposition]{symPower}
    For $m\leq 4$ we have
    \[
        dim\Ext^p(\Sym^m\mQ,\Sym^m\mQ)=
        \begin{cases}
        1&p=0,4\\
        0&p=1,3\\
        3\left(\frac{3m^2+12m-20}{20}\right)^2r^2_m-2&p=2.
        \end{cases}
    \]
More precisely, for $0\leq m\leq 4$,\ $\Ext^2(\Sym^m\mQ,\Sym^m\mQ)$ is the direct sum of vector spaces in the columns of \cref{tab:ext2Sym}.
\begin{table} 
    \centering
    \begin{adjustbox}{max width=\textwidth}
    \begin{tabular}{|c|c|c|c|c|}
        \hline
         $m=0$ & $m=1$ & $m=2$ & $m=3$ & $m=4$ \\ 
        \hline
         $\mathbb{C}$ & $\mathbb{C}$ & $\mathbb{C}$ & $\mathbb{C}$ & $\mathbb{C}$\\ 
        \hline
         - & - & $\Schur_{1,1,1,1,1,1,1,1,0,0}(V_{10}\oplus V_{10}^\vee)$ & $\Schur_{1,1,1,1,1,1,1,1,0,0}(V_{10}\oplus V_{10}^\vee)$ & $\Schur_{1,1,1,1,1,1,1,1,0,0}(V_{10}\oplus V_{10}^\vee)$ \\ 
        \hline
         - & - & $\Schur_{2,1,1,1,1,1,1,1,1,0}V_{10}$ & $\Schur_{2,1,1,1,1,1,1,1,1,0}V_{10}$ & $\Schur_{2,1,1,1,1,1,1,1,1,0}V_{10}$ \\ 
        \hline
         - & - & - & $\Sigma_{3,0,0,0,0,0,0,0,0,0}(V_{10}\oplus V_{10}^\vee)$ & $\Sigma_{3,0,0,0,0,0,0,0,0,0}(V_{10}\oplus V_{10}^\vee)$ \\ 
        \hline
         - & - & - & $\Sigma_{2,1,1,1,1,1,1,0,0,0}(V_{10}\oplus V_{10}^\vee)$ & $\Sigma_{2,1,1,1,1,1,1,0,0,0}(V_{10}\oplus V_{10}^\vee)$ \\ 
        \hline
        - & - & - & $V_{10}^{\oplus 2}$ & $V_{10}^{\oplus 2}$ \\ 
        \hline
        - & - & - & $\Sigma_{2,2,2,2,2,2,2,2,1,0}(V_{10}\oplus V_{10}^\vee)$ & $\Sigma_{2,2,2,2,2,2,2,2,1,0}(V_{10}\oplus V_{10}^\vee)$ \\ 
        \hline
        - & - & - & $\Sigma_{2,2,1,1,1,1,1,1,0,0}V_{10}$ & $\Sigma_{2,2,1,1,1,1,1,1,0,0}V_{10}$ \\ 
        \hline
        - & - & - & -& $\Schur_{2,1,1,0,0,0,0,0,0,0}(V_{10}\oplus V_{10}^\vee)$\\
        \hline
        - & - & - & -& $\Schur_{3,2,2,2,2,2,2,2,2,0}(V_{10}\oplus V_{10}^\vee)$\\
        \hline
        - & - & - & -& $\Schur_{2,2,2,2,1,1,0,0,0,0}V_{10}$ \\
        \hline
    \end{tabular}
    \end{adjustbox}
    \newline 
    \caption{}
    \label[table]{tab:ext2Sym}
    \end{table}
\end{proposition}
\begin{proof}

Here we produced the vector space explicitly using Borel--Weil--Bott’s theorem. In particular:
\begin{itemize}
    \item $(m=0)$ By looking at the table \cref{tab:sym-transposed}, we have that the partitions associated to non-acyclic bundles for $m=0$ are $(0,0,0,0|0,0,0,0,0,0)$ and  $(0,0,0,0|7,7,7,3,3,3)$, which corresponds to $\of_{\Gr(6,10)}$ and to a factor of $\bigwedge^{10}\bigwedge^3\widetilde{\mU}$. We recover $\Ext^0(\of_X,\of_X)=H^0(X,\of_{\Gr(6,10)})=\mathbb{C}$, and by symmetry also $\Ext^4(\of_X,\of_X)$. Moreover $\Ext^2(\of_X,\of_X)=H^{12}(\Gr(6,10),\Schur_{0,0,0,0|7,7,7,3,3,3})=\mathbb{C}$.\\
    \item $(m=1)$ All the partitions in the table are associated with acyclic bundles. By \cref{rmk:sym-skip} we get that $\Ext^p(\mQ,\mQ)=\Ext^p(\of_X,\of_X)$. This is also a case already known.\\
    \item $(m=2)$ The partitions in \cref{tab:sym-transposed} associated to a non-acyclic bundle for $m=2$ are $(4,2,2,0|7,7,5,5,3,3)$, factor of $E_2^6$, and $(4,2,2,0|10,8,8,6,6,4)$, factor of $E_2^{10}$. They only contribute to the $\Ext$-groups in the $\Ext^2$. In particular: 
    
    \begin{align*}
    \Ext^0(\Sym^2\mQ,\Sym^2\mQ)~&=\Ext^4(\Sym^2\mQ,\Sym^2\mQ)=\Ext^0(\mQ,\mQ)=\mathbb{C},\\
    \Ext^1(\Sym^2\mQ,\Sym^2\mQ)~&=\Ext^3(\Sym^2\mQ,\Sym^2\mQ)=0,\\  
    \Ext^2(\Sym^2\mQ,\Sym^2\mQ)~&=H^{8}(\Gr(6,10),\Schur_{4,2,2,0|7,7,5,5,3,3})^{\oplus 2}\oplus H^{12}(\Gr(6,10),\Schur_{4,2,2,0|10,8,8,6,6,4})\oplus\Ext^2(\mQ,\mQ)=\\&=(\Schur_{1,1,1,1,1,1,1,1,0,0}V_{10})^{\oplus 2}\oplus\Schur_{2,1,1,1,1,1,1,1,1,0}V_{10}\oplus\Ext^2(\mQ,\mQ).
    \end{align*} 
    The second copy $\Schur_{1,1,1,1,1,1,1,1,0,0}V_{10}^\vee$ comes from the symmetric of $\Schur_{4,2,2,0|7,7,5,5,3,3}$ in $\bigwedge^{16}\bigwedge^3\widetilde{\mU}$.\\
    \item $(m=3)$ In this case we have five partitions in \cref{tab:sym-transposed} associated with non-acyclic bundles: 
    \begin{align*}
        (6,3,3,0|4,4,4,3,3,3) \text{ in } E_3^1,\quad (6,3,3,0|8,8,6,6,4,4) &\text{ in }E_3^6, \quad(6,3,3,0|8,8,6,6,6,5) \text { in }E_3^7,
    \end{align*}
    \vspace{-2.5em}
    \begin{align*}(6,3,3,0|9,9,9,7,7,4) \text{ in}\ E_3^9, \quad&\text{and}\quad
        (6,3,3,0|11,9,9,7,7,5) \text{ in } E_3^{10}. 
    \end{align*}
    Using Borel--Weil--Bott’s theorem, we note that they contribute only to the $\Ext^2$-group. In particular,

    \begin{align*}
    \Ext^0(\Sym^3\mQ,\Sym^3\mQ)~&=\Ext^4(\Sym^3\mQ,\Sym^3\mQ)=\Ext^0(\Sym^2\mQ,\Sym^2\mQ)=\mathbb{C},\\ \Ext^1(\Sym^3\mQ,\Sym^3\mQ)~&=\Ext^3(\Sym^3\mQ,\Sym^3\mQ)=\Ext^1(\Sym^2\mQ,\Sym^2\mQ)=0,\\
    \Ext^2(\Sym^3\mQ,\Sym^3\mQ)~&=H^{3}(\Gr(6,10),\Schur_{6,3,3,0|4,4,4,3,3,3})^{\oplus 2}\oplus H^{8}(\Gr(6,10),\Schur_{6,3,3,0|8,8,6,6,4,4})^{\oplus 2}\oplus\\ 
    &\oplus H^{9}(\Gr(6,10),\Schur_{6,3,3,0|8,8,6,6,6,5})^{\oplus 2}\oplus H^{11}(\Gr(6,10),\Schur_{6,3,3,0|9,9,9,7,7,4})^{\oplus 2}\oplus\\ 
    &\oplus H^{12}(\Gr(6,10),\Schur_{6,3,3,0|11,9,9,7,7,5})\oplus\Ext^2(\Sym^2\mQ,\Sym^2\mQ)=\\ 
    &=(\Sym^3 V_{10})^{\oplus 2}\oplus(\Schur_{2,1,1,1,1,1,1,0,0,0} V_{10})^{\oplus 2}\oplus V_{10}^{\oplus 2}\oplus(\Schur_{2,2,2,2,2,2,2,2,1,0} V_{10})^{\oplus 2}\oplus\\ 
    &\oplus\Schur_{2,2,1,1,1,1,1,1,0,0} V_{10}\oplus\Ext^2(\Sym^2\mQ,\Sym^2\mQ).
    \end{align*}
    The double copies of cohomology space come from the symmetric component in the Koszul Complex as explained in \cref{rmk:dWtoSchur}.
    \item $(m=4)$ We have three partitions associated to non-acyclic bundles: $(8,4,4,0|10,10,7,7,7,7)$ in $E_4^8$, $(8,4,4,0|10,10,10,8,8,5)$ in $E_4^9$, and $(8,4,4,0|11,11,11,7,7,7)$ in $E_4^{10}$. This again, by Borel--Weil--Bott’s theorem, gives:

    \begin{align*}
    \Ext^0(\Sym^4\mQ,\Sym^4\mQ)~&=\Ext^4(\Sym^4\mQ,\Sym^4\mQ)=\Ext^0(\Sym^3\mQ,\Sym^3\mQ)=\mathbb{C},\\ \Ext^1(\Sym^4\mQ,\Sym^4\mQ)~&=\Ext^3(\Sym^4\mQ,\Sym^4\mQ)=\Ext^1(\Sym^2\mQ,\Sym^2\mQ)=0,\\
    \Ext^2(\Sym^4\mQ,\Sym^4\mQ)~&= H^{10}(\Gr(6,10),\Schur_{8,4,4,0|10,10,7,7,7,7})^{\oplus 2}\oplus H^{11}(\Gr(6,10),\Schur_{8,4,4,0|10,10,10,8,8,5})^{\oplus 2}\oplus \\
    &\oplus H^{12}(\Gr(6,10),\Schur_{8,4,4,0|11,11,11,7,7,7})\oplus\Ext^2(\Sym^3\mQ,\Sym^3\mQ)=\\
    &=(\Schur_{2,1,1,0,0,0,0,0,0,0}V_{10})^{\oplus 2}\oplus (\Schur_{3,2,2,2,2,2,2,2,2,0}V_{10})^{\oplus 2}\oplus\\ &\oplus\Schur_{2,2,2,2,1,1,0,0,0,0}V_{10}\oplus\Ext^2(\Sym^3\mQ,\Sym^3\mQ) 
    \end{align*}
We obtain the dimensions of $\Ext$-spaces by computing the dimensions of the vector spaces we found. In this way, we could describe explicitly the $\Ext$-spaces, but we could also use \cref{ChernSym}. If we restrict ourselves to $0\leq m\leq 4$, then we can note that we do not have any contribution for $$\Ext^1(\Sym^m\mQ,\Sym^m\mQ)=\Ext^3(\Sym^m\mQ,\Sym^m\mQ),$$ hence it is always $0$ in these cases. Moreover the only contribution to $\Ext^0(\Sym^m\mQ,\Sym^m\mQ)$ and $\Ext^4(\Sym^m\mQ,\Sym^m\mQ)$ is from $H^0(X,\of_X)=\mathbb{C}$. Hence, by \cref{ChernSym}: \[
\mathrm{dim}\Ext^2(\Sym^m\mQ,\Sym^m\mQ)=\chi(\Sym^m\mQ,\Sym^m\mQ)-2=3\left(\frac{3m^2+12m-20}{20}\right)^2r^2_m-2.
\]

\end{itemize}

\end{proof}

\begin{rmk}
    Note that for $m\geq 7$ all the partitions in \cref{tab:sym-transposed} are associated with non-acyclic bundles. In particular, one can try to mimic the method used in \cite[Prop. 2.3]{FatOno}, to obtain a stronger result. The problem is that from $m\geq 5$, $\mathrm{Sym}^m\mQ$ is indeterminate because of the $E_5^q$'s, as we show in the following example. Thus, it is not possible to conclude, even if it is reasonable that \cref{symPower} holds even for $m>4$.
\end{rmk}

\begin{ex}\textbf{$\Ext(\Sym^5\mQ,\Sym^5\mQ)$.} We proceed as before. First note:
\[
\Ext^p(\Sym^5\mQ,\Sym^5\mQ)=H^p(X,E_5)\oplus\Ext^p(\Sym^4\mQ,\Sym^4\mQ)
\]
So we study, as before, only the cohomology of $E_5$. We recover the Koszul complex:
\[
0\to E_5^{20}\to...\to E_5^{2}\to E_5^1\to E_5^0\to E_5\to 0. 
\]
Let $0\leq q\leq 10$, then we have nine partitions in \cref{tab:sym-transposed} that are associated to non-acyclic factors of the $E_5^q$'s, we recollect them in the following \cref{tab:symFive}:

\begin{table} 
    \centering
    \begin{adjustbox}{max width=\textwidth}
    \begin{tabular}{|c|c|c|c|c|}
        \hline
         q=6 & - & $(10,5,5,0|10,10,7,7,7,7)$ & - & - \\ 
        \hline
         q=7 & $(10,5,5,0|11,10,8,8,7,7)$ & $(10,5,5,0|10,10,8,8,8,7)$ & - & - \\ 
        \hline
         q=8 & $(10,5,5,0|11,11,8,8,8,8)$ & - & - & - \\ 
        \hline
         q=9 & $(10,5,5,0|12,11,11,8,8,7)$ & - & - & - \\ 
        \hline
         q=10 & $(10,5,5,0|13,11,11,9,9,7)$ & $(10,5,5,0|13,11,11,9,8,8)$ & $(10,5,5,0|12,12,12,8,8,8)$ & $(10,5,5,0|12,12,11,9,9,7)$ \\ 
        \hline
    \end{tabular}
    \end{adjustbox}
    \newline 
    \caption{}
    \label[table]{tab:symFive}
    \end{table}

By Borel--Weil--Bott’s theorem and \cref{rmk:dWtoSchur} we obtain the following non-trivial cohomologies:
\begin{align*}
 H^{10}(\Gr(6,10),E^6_5)~&=H^{14}(\Gr(6,10),E^{14}_6)=\Schur_{4,2,1,1,1,0,0,0,0,0}V_{10}\\
 H^{10}(\Gr(6,10),E^{7}_5)~&=H^{14}(\Gr(6,10),E^{13}_5)=\Schur_{4,2,1,1,1,1,0,0,0,0}V_{10}\oplus\Schur_{4,1,1,1,1,1,1,1,0,0}V_{10}\\
 H^{10}(\Gr(6,10),E^{8}_5)~&=H^{14}(\Gr(6,10),E^{12}_5)=\Schur_{4,2,2,1,1,1,1,1,1,0}V_{10}\\
 H^{12}(\Gr(6,10),E^{9}_5)~&=H^{12}(\Gr(6,10),E^{11}_5)=\Schur_{4,3,2,2,2,2,1,1,0,0}V_{10}\\
 H^{12}(\Gr(6,10),E^{10}_5)~&=\Schur_{4,4,2,2,2,2,2,2,0,0}V_{10}\oplus\Schur_{4,4,2,2,2,2,2,1,1,0}V_{10}\oplus\Schur_{4,3,3,3,2,2,1,1,1,0}V_{10}\oplus\Schur_{4,3,3,2,2,2,2,2,0,0}V_{10}
\end{align*}

Since we have non-trivial cohomologies on the same degree, the short exact sequences that arise in cohomology are not trivially solved. This is the reason we can not conclude the argument for $m>4$.
\end{ex}

\section{Comparison with atomicity}\label[section]{sec:atom}
Another notion introduced to study the modularity properties of sheaves on \hk manifolds was given by Beckmann in \cite{Bec25}. There, a coherent sheaf on a \hk manifold $X$ is called \textit{atomic} if its Mukai vector "behaves" as a vector in the extended Mukai lattice of $X$ with respect to the action of the LLV algebra $\mathfrak{g}(X)$, see Definition 1.1, \textit{op. cit.}. We provide a definition below which, in the $\textrm{K3}^{[2]}$-type case, is equivalent to the general one.

In this section, we will prove that the only atomic bundles among the Schur functors of $\mQ$ are given by the symmetric powers $\mathrm{Sym}^m(\mQ)$. Note that, in general, any atomic torsion-free sheaf is modular by \cite[Proposition 1.5]{Bec25}. Our main result is the following:
\begin{proposition}\label[proposition]{atomicity}
    Let $X\subseteq \mathrm{Gr}(6,10)$ be a very general Debarre--Voisin \hk manifold, and let $\lambda = (m,t,s,0)$ with $m\geq t+s$, then $\Sigma_\lambda \mQ$ is atomic if and only if $\lambda=(m,0,0,0)$.
\end{proposition}

We divide the proof into steps.
\subsection{Proof of \texorpdfstring{\cref{atomicity}}{atomicity}}\label[section]{proof7.1}
We recall some definitions and notations from \cite{Bec25}. Note that the following statements are specializations of more general ones to the $\textrm{K3}^{[2]}$-type case.
\begin{definition} Let $X$ a \hk manifold of $\textrm{K3}^{[2]}$-type.
\begin{itemize}
    \item The rational \textit{extended Mukai Lattice} is the graded vector space 
    \[\widetilde{H}(X,\mathbb{Q}) := \mathbb{Q}\alpha \oplus H^2(X,\mathbb{Q})\oplus \mathbb{Q}\beta\]
    endowed with a bilinear form $\widetilde{q}$ which restricts to the Beauville--Bogomolov--Fujiki form $q_X$ on $H^2(X,\mathbb{Q})$ and such that $\widetilde{q}(\alpha,\alpha) = \widetilde{q}(\beta,\beta) = \widetilde{q}(\alpha,H^2(X,\mathbb{Q})) = \widetilde{q}(\beta,H^2(X,\mathbb{Q})) = 0$ and $\widetilde{q}(\alpha,\beta)=-1$.
    \item There is an isometric embedding $\psi:H^*(X,\mathbb{Q}) \hookrightarrow \mathrm{Sym}^2\widetilde{H}(X,\mathbb{Q})$. We denote by $T:\mathrm{Sym}^2\widetilde{H}(X,\mathbb{Q}) \longrightarrow H^*(X,\mathbb{Q})$ the correspondent orthogonal projection. 
    \item A coherent sheaf $\mathcal{E}$ on $X$ is \textit{atomic} if there exists a vector $\widetilde{v}\in \widetilde{H}(X,\mathbb{Q})$ such that the Mukai vector of $\mathcal{E}$, $v(\mathcal{E}):= \ch(\mathcal{E})\sqrt{\mathrm{td}_X}$, is a rational multiple of $T(\widetilde{v}\cdot\widetilde{v})$. In this case, $\widetilde{v}$ is called \textit{extended Mukai vector} for $\mathcal{E}$.
\end{itemize}
    
\end{definition}

\begin{lemma}
    Following the notations of \cref{AboutSym}, let 
    \[
    \widetilde{v}_m := \left(r_m,\frac{m}{4}r_m\ch_1(\mQ),\frac{2m^2-3m+5}{4}r_m\right)\in \widetilde{H}(X,\mathbb{Q}),
    \]
    then $\widetilde{v}_m $ is an extended Mukai vector for $\mathrm{Sym}^m(\mQ)$. In particular,  $\mathrm{Sym}^m(\mQ)$ is atomic for every $m\geq 1$.
\end{lemma}
\begin{proof}
    We first compute the Mukai vector  $v(\mathrm{Sym}^m(\mQ)):= \ch(\mathrm{Sym}^m(\mQ))\sqrt{\mathrm{td}_X}$, using item (2) of \cref{ChernSym}, item (3) of \cref{toddIHS} and \cref{ChernRelations}.
    \begin{align*}
    v(\mathrm{Sym}^m(\mQ))= &\left(r_m,\frac{m}{4}r_m\ch_1(\mQ),\right.\\
    &\frac{m^2-m}{40}r_m \ch_1(\mQ)^2 + \frac{m^2+4m}{20}r_m\ch_2(\mQ)+\frac{1}{24}r_m\ch_1(\mQ)^2-\frac{1}{3}r_m \ch_2(\mQ),\\
     &-\frac{2m^3-3m^2}{4}r_m \ch_3(\mQ)-\frac{11m}{4}r_m \ch_3(\mQ) + \frac{3}{2}r_m \ch_3(\mQ),\\
     &-\frac{10m^4-30m^3+21m^2-6m}{20} r_m \ch_4(\mQ) -\frac{121}{20}(m^2-m)r_m \ch_4(\mQ)+\\
     &-\frac{33}{40}(m^2+4m)r_m \ch_4(\mQ) +\frac{33}{10}(m^2-m)r_m \ch_4(\mQ) + (m^2+4m)r_m \ch_4(\mQ)+ \\
     &\left.-\frac{25}{8}r_m \ch_4(\mQ)\right).
    \end{align*}
    If we denote by $\ch_1(\mQ)^\vee \in H^6(X,\mathbb{Q})$ the dual class as in the proof of \cref{ChernRelations}, it follows that 
    \begin{align*}
    v(\mathrm{Sym}^m(\mQ))= &\left(r_m,\frac{m}{4}r_m\ch_1(\mQ),\right.\\
    &\frac{3m^2-3m+5}{120}r_m \ch_1(\mQ)^2 + \frac{3m^2+12m-20}{60}r_m\ch_2(\mQ),\\
     &\frac{2m^3-3m^2+5m}{16}r_m \ch_1(\mQ)^\vee,\\
     &\left.\frac{4m^4-12m^3+29m^2-30m+25}{32}r_m \right).
    \end{align*}
    On the other hand, we can compute the image of $\widetilde{v}_m\cdot\widetilde{v}_m$ using the following formula (\cite[Lemma 1.19]{FatOno})
    \[
    T((r,\ell,s)\cdot (r,\ell,s)) = \left(r,\ell,\frac{1}{2r}\left(\ell^2-\frac{\widetilde{q}((r,\ell,s)}{30} \mathsf{c}_2(X)\right),\frac{s}{r}\ell^\vee,\frac{s^2}{2r}\right).
    \]
    For $r=r_m$, $\ell = \ch_1(\mathrm{Sym}^m(\mQ))$ and $s=\frac{2m^2-3m+5}{4}r_m$, the formula yields
    \begin{align*}
    T(\widetilde{v}_m\cdot\widetilde{v}_m) = &\left(r_m,\frac{m}{4}r_m\ch_1(\mQ),\frac{1}{2r_m}\left(\frac{m^2}{12}r_m^2\ch_1(\mQ)^2-\frac{\widetilde{q}(\widetilde{v}_m)}{30} \mathsf{c}_2(X)\right),\right.\\
    &\left.\frac{2m^3-3m^2+5m}{16}r_m\ch_1(\mQ)^\vee,\frac{4m^4-12m^3+29m^2-30m+25}{32}r_m\right).
    \end{align*}
    Moreover, $\widetilde{q}(\widetilde{v}_m)= \frac{11m^2}{8}r_m^2-2\frac{2m^2-3m+5}{4}r_m^2=\frac{3m^2+12m-20}{8}r_m^2$. After a short computation, one obtains the desired 
    \[ T(\widetilde{v}_m\cdot\widetilde{v}_m)=v(\mathrm{Sym}^m(\mQ)).
    \]
    Thus, $\widetilde{v}_m$ is an extended Mukai vector for $\mathrm{Sym}^m(\mQ)$ and thus $\mathrm{Sym}^m(\mQ)$ is atomic.
\end{proof}

Conversely, suppose that $\Sigma_\lambda \mQ$ is atomic and let $\widetilde{v}_m$ be a extended Mukai vector for $\Sigma_\lambda \mQ$. Then, by \cite[Proposition 3.12 and Theorem 3.17]{Bot24}, we have 
\begin{align*}
    \Delta(\Sigma_\lambda \mQ) &= \frac{\widetilde{q}(\widetilde{v}_m)+\frac{5}{2}r^2}{30}\mathsf{c}_2(X),\\
  \chi(\Sigma_\lambda \mQ,\Sigma_\lambda \mQ)&= 3\left(\frac{\widetilde{q}(\widetilde{v}_m)}{\frac{5}{2}r^2}\right)^2 r^2.
\end{align*}

On the other hand, by \cref{SchurIsModular} and \cref{SchurIsModular2}, we have 
\begin{align*}
    \Delta(\Sigma_\lambda \mQ)&= \frac{\delta}{4} r^2\mathsf{c}_2(X),\\
    \chi(\Sigma_\lambda \mQ,\Sigma_\lambda \mQ)&= 3\left(1+\frac{-276\delta-1936\ell^4+1320\delta\ell^2+207\delta^2-8\xi}{48}\right)r^2.
\end{align*}

In particular, $P(m,t,s) := 1+\frac{-276\delta-1936\ell^4+1320\delta\ell^2+207\delta^2-8\xi}{48}$ has to be equal to $(3\delta -1)^2$. In analogy to \cite[Proposition 6.4]{FatOno}, a direct computation shows that, up to duality and tensor product with line bundles, this forces the partition $\lambda=(m,t,s,0)$ to be $(m,0,0,0)$. 

The proof of \cref{atomicity} is therefore complete.

\frenchspacing

\newcommand{\etalchar}[1]{$^{#1}$}

\end{document}